\documentclass[12pt]{amsart}
\usepackage{amsmath,amssymb,amsthm,geometry}
\usepackage{amsrefs}

    \usepackage[pdftex]{graphicx}
    \usepackage[colorlinks=true, pdfstartview=FitH, linkcolor=blue,
        citecolor=blue, urlcolor=blue]{hyperref}

    \newfont{\footsc}{cmcsc10 at 8truept}
    \newfont{\footbf}{cmbx10 at 8truept}
    \newfont{\footrm}{cmr10 at 10truept}
    \title[Comparative prime number theory: a survey]{Comparative Prime Number Theory:\\ 
    A survey}

\author{Greg Martin}
\author{Justin Scarfy}
\address{University of British Columbia\\ Department of Mathematics \\ Room 121\\ 1984 Mathematics Road\\Vancouver, BC Canada V6T 1Z2}
\email{gerg@math.ubc.ca} 
\email{scarfy@ugrad.math.ubc.ca}

\thanks{The second of us, an undergraduate student at the University of British Columbia, is very grateful for Professor Greg Martin's  introduction to this fascinating topic, and for his guidance and continuous encouragement throughout this project.}


    \date{\small 
 February 14, 2012}

\newcommand{\h}{\frac 1 2}
\renewcommand{\mod}[1]{{\ifmmode\text{\rm\ (mod~$#1$)}\else\discretionary{}{}{\hbox{ }}\rm(mod~$#1$)\fi}}
\newcommand{\Li}{\rm Li}
\newcommand{\Var}{\rm Var}

\newcommand{\li}{\rm li}
\newcommand{\ord}{\rm ord}
\newcommand{\N}{{\mathbb N}}

\newcommand{\Q}{{\mathbb Q}}
\newcommand{\R}{{\mathbb R}}

\newcommand{\A}{{\mathcal A}}

\numberwithin{equation}{section}

\newtheorem{thm}{Theorem}[subsection]

\newtheorem{cnj}[thm]{Conjecture}
\newtheorem{dfn}[thm]{Definition}
\newtheorem{rmk}[thm]{Remark}
\newtheorem{xmp}[thm]{Example}
\newtheorem{pbm}[thm]{Problem}

\numberwithin{thm}{section} 

\begin{document}
\begin{abstract}
Comparative prime number theory is the study of the {\em{discrepancies}} of distributions when we compare the number of primes in different residue classes. This work presents a list of the problems being investigated in comparative prime number theory, their generalizations, and an extensive list of references on both historical and current progresses. 
\end{abstract}
\maketitle
\section{Introduction}\label{intro section}
In a letter between P. Chebyshev and M. Fuss, dated 1853~\cite{1853.Chebyshev}, the former indicates (without proof):
For a positive continuous decreasing function $f$, the series
\begin{equation} \label{eq.cheb} 
\sum_{p \; {\rm {odd\;prime}}} (-1)^{\frac {(p+1)}{2}} f(p): = f(3)-f(5)+f(7)+f(11)-f(13)-f(17)+\hdots
\end{equation}
diverges. In particular, when $f(x)=e^{-10x}$, the series \eqref{eq.cheb} tends to infinity.
The significance of this assertion is to say that there {\em{should be}} more primes in the residue class 3 modulo 4 than in the residue class 1 modulo 4, despite the fact that Dirichlet in 1837 proved that for any $a, k$ with $(a, k)=1$ there are infinitely many primes $p$ with $p\equiv a \mod k$.
Hardy, Littlewood, and Landau in 1918 proved that Chebyshev's assertion is equivalent to the problem of whether the function 
\begin{equation*}
L(s):= \sum_{n=0}^{\infty} \frac{(-1)^n}{(2n+1)^s} \; \; \; \; \; (s=\sigma + it)
\end{equation*}
vanishes or not in the half-plane $\sigma > \h$. (The necessity was shown by Landau~\cite{1918.Landau.1} and the sufficiency by Hardy--Littlewood~\cite{1918.Littlewood}, with a simpler proof by Landau~\cite{1918.Landau.2}.)

However, Littlewood~\cite{1918.Littlewood} in 1914 showed that the number of primes in the residue class 3 modulo 4 and the number of primes in the residue class 1 modulo 4 ``race'', taking turns to be in the lead.  
On the other hand, the number of primes in the residue class 1 seems to take the lead in the race only a ``negligible'' amount of time, and this phenomenon is known as {\bf Chebyshev's bias}.  
To illustrate precisely what Littewood had proven and further developments on this topic, we need the aid of the following notations:

Throughout this paper, $p$ will always be an {\bf{odd} prime}. As usual, for a positive integer $k$ with $(k, l)=1$:
\begin{dfn}
\[ 
\pi(x; k, l) := \sum_{\substack {p\leq x \\ p\equiv l \mod k}} 1
\]
\end{dfn}
We note that $\pi(x)=\pi(x; 1, 1)$, the prime counting function up to $x$.

By the prime number theorem for arithmetic progressions, the functions $\pi(x; k, l)$ with fixed $k$ and $(l, k)=1$ are all asymptotically $x/(\varphi(k)\log x)$, where $\varphi(k)$ is the Euler's totient function. As Chebyshev investigated, the difference between the functions $\pi(x; k, l)$ for fixed $k$ exhibits interesting behaviours:
\begin{dfn}
\[
\delta_\pi(x; k, l_1, l_2):=\pi(x; k, l_1)-\pi(x; k, l_2)
\]
\end{dfn}
\begin{xmp}
In Chebyshev's case, $\delta_\pi(x; 4, 3, 1)$ is negative for the first time when $x=26$,$861$~\cite{1957.Leech}, and $\delta_\pi(x; 3, 2, 1)$ is negative for the first time at $x=608$,$981$,$813$,$029$~\cite{1977.Bays}.
\end{xmp}

\begin{xmp}
Littlewood ~\cite{1918.Littlewood} proved in 1914 that, in the above notations, $\delta_\pi(x; 4, 3, 1)$ and $\delta_\pi(x; 3, 2, 1)$ switch their signs infinitely many times.
\end{xmp}

\begin{xmp}
Phragm\'en ~\cite{1891.Phragmen} proved the existence of an unbounded sequence $x_1 < x_2 < x_3 < \cdots$ such that
\[
\frac{\pi(x_\nu; 4, 3)- \pi(x_\nu; 4, 1)}{{\sqrt {x_\nu}}/{\log x_\nu}} \rightarrow 1
\]
\end{xmp}

In their series of papers published between 1962 and 1972, S. Knapowski and P. Tur\'an ~\cite{1962.Knaposwki_1,1964.Turan_1} list a number of problems that generalized Littlewood's theorem and also attempted to compare $\pi(x; k, l_1)$ and $\pi(x, k, l_2)$, with the assumption that $l_1 \not\equiv l_2 \mod k$ and $(l_1, k)=(l_2, k)=1$:

\begin{pbm}[Infinity of sign changes]\label{P:S-C}
For which $(l_1, l_2)$-pairs  does the function $\delta_\pi(x; k, l_1, l_2)$ change its sign infinitely often?
\end{pbm}

\begin{pbm}[Big sign changes]\label{P:B-S}
Given $\epsilon >0$, do there exist two sequences 
\begin{align*}
x_1&<x_2<x_3<\cdots \rightarrow  \infty\\
y_1&<y_2<y_3<\cdots \rightarrow  \infty
\end{align*}
such that
\begin{align*}
\pi(x_\nu ; k, l_1)-\pi(x_\nu ; k,  l_2) &> x_\nu ^ {1/2 - \epsilon} \\
\pi(y_\nu ; k, l_1)-\pi(y_\nu ; k,  l_2) &< -y_\nu ^ {1/2 - \epsilon}?
\end{align*}
The use of the function $x^{1/2 -\epsilon}$ is motivated by the fact that if GRH for $k$ (see Conjecture~\ref{GRH} below) is true, then the inequality
\[
|\delta_\pi(x; k, l_1, l_2)|=O(x^{1/2}\log x)
\]
holds for $x\ge2$.
\end{pbm}


\begin{pbm}[Localized sign changes]\label{P:L-S}
Prove that there exists an $T>T_0(k)$ and suitable $A(T)<T$ such that the function $\delta_\pi(x; k, l_1, l_2)$ changes sign in the interval 
$$
A(T)\leq x \leq T.
$$
\end{pbm}

\begin{pbm}[Localized big sign changes]\label{P:L-B}
Prove that for $T> T_0(k)$ and suitable $A(T)<T$, the functions $\delta_\pi(x; k, l_1, l_2)$ satisfy both the inequalities
\begin{align*}
\max_{A(T)\leq x \leq T}\delta_\pi(x; k, l_1, l_2) &>  \frac{T^{1/2}}{\Phi(T)}\\
\min_{A(T)\leq x \leq T}\delta_\pi(x; k, l_1, l_2) &< {-} \frac{T^{1/2}}{\Phi(T)},
\end{align*}
where $\Phi(x)$ is a positive function satisfying 
$$\lim_{x\to \infty}\frac{\log\Phi(x)}{\log x}= 0.$$
\end{pbm}

\begin{pbm}[First sign change]\label{P:F-S}
For what function $a(k)$ can we assert that for each $(l_1, l_2)$-pair with $l_1 \neq l_2$, all functions in
$
\delta_\pi(x; k, l_1, l_2)$
vanish at least once in
$
1 \leq x \leq a(k)?
$
\end{pbm}

\begin{pbm}[Asymptotic behaviour of sign changes]\label{P:A-B}
Let  $w_\pi(T; l_1, l_2)$ denote the number of sign changes of $\delta_\pi(x; k, l_1, l_2)$ in the interval $[1, T]$. What is the asymptotic behaviour of $w_\pi(T; l_1, l_2)$ as $T \rightarrow \infty$?
\end{pbm}


\begin{pbm}[Race-problem of Shanks--R\'enyi]\label{P:P-R}
For each permutation $\{l_1, l_2, l_3, \hdots, l_{\varphi(k)}\}$ of the set of reduced residue classes modulo $k$, do there exist infinitely many integers $m$ with
\[
\pi(m; k, l_1) < \pi(m; k, l_2) < \pi(m; k, l_3) < \cdots  < \pi(m; k, l_{\varphi(k)})?
\]
G.\ G.\ Lorentz noticed the fact that comparison of primes of any two arithmetical progressions mod $k_1$ and $k_2$ $(k_1 \neq k_2)$ is not trivial in the case when
\begin{equation*}
\varphi(k_1)= \varphi(k_2)
\end{equation*}
and analogous problems occur for moduli $k_1, k_2, k_3, \hdots, k_r$ with 
\begin{equation*}
\varphi(k_1)= \varphi(k_2)= \cdots = \varphi(k_r),
\end{equation*}
\end{pbm}

\begin{dfn}
Define
\[
\Li(x):=\int_0^x\frac{dt}{\log t}.
\]
\end{dfn}

\begin{pbm}[Littlewood generalizations]\label{P:L-G}
Do there exist infinitely many integers $m_\nu$ such that for $j=1, 2, 3, \hdots \varphi(k)$, we simultaneously have
$$
\pi(m_\nu, k, l_j)> \frac {\Li(m_\nu)}{\varphi(k)} ?
$$
and if the assertion is valid, what are the distribution-properties of the sequence~$m_\nu$?
\end{pbm}

\begin{pbm}[Average preponderance problems]\label{P:A-P}
Denote by $N_\pi(x)$ the number of integers $n\leq x$ with the property $\delta_\pi(n; 4, 3, 1)>0$.
Does the relation
\begin{equation*}
\lim_{x \to \infty} \frac{N_\pi(x)}x = 0
\end{equation*}
hold? In other words, does the set of integers $n$ with the property $\delta_\pi(n; 4, 3, 1)>0$ have density~$0$?
\end{pbm}

In the previous problems, the number of \emph{all} primes $\leq x$ in a fixed progression occurred. One can imagine that one can much better locate relatively small intervals where the primes of some progression preponderate.
\begin{pbm}[Strongly localized accumulation problems]\label{P:S-L}
When $T$ is sufficiently large, is it true that for suitable $T \leq U_1 < U_2 \leq 2T$, we have
\begin{equation*}
\sum_{\substack{U_1\leq p \leq U_2 \\ p \equiv l_1 \mod k}}1 - \sum_{\substack{U_1\leq p \leq U_2 \\ p \equiv l_2 \mod k}}1 > \frac{\sqrt T}{\Phi(T)}?
\end{equation*}
\end{pbm}
where $\Phi(x)$ shares the same property as in Problem~\ref{P:L-B}

\begin{pbm}[Union problem]\label{P:U-P}
For a given modulus $k$, do there exist two disjoint subsets $A$ and $B$, consisting of the same number of reduced residue classes, such that 
\begin{equation*}
\sum_{\substack{p\in A \\ p\leq x}} 1 \geq \sum_{\substack{p\in B \\ p\leq x}} 1
\end{equation*}
for all sufficiently large $x$?
\end{pbm}

\begin{rmk}\label{l_1 and l_2}
One can expect that there are ``more'' primes in the residue class $l_1\mod k$ than $l_2 \mod k$ if and only if the number of incongruent solutions of the congruence
\begin{equation}\label{l_1}
x^2 \equiv l_1 \mod k
\end{equation} 
is less than that of the congruence
\begin{equation}\label{l_2}
x^2 \equiv l_2 \mod k
\end{equation}
\end{rmk} 

Besides the functions $\pi(x; k, l)$, the distributions of primes in arithmetic progressions can be studied by some other functions that are easier to work with. Let $\Lambda(n)$ denote the von Mangoldt Lambda function, namely:
\begin{dfn}
\[
\Lambda(n):= 
\begin{cases}  
\log{p} & {\rm if }\;\; n=p^k\\
      0 & {\rm otherwise} 
\end{cases}
\]
\end{dfn}
And thus the following functions are studied:
\begin{dfn}
\begin{align*}
\psi(x; k, l) :=& \sum_{\substack {n\leq x \\ n\equiv l \mod k}} \Lambda(n) \\
\vartheta(x; k, l) :=& \sum_{\substack {p\leq x \\ p\equiv l \mod k}} {\log p}  \\
\Pi(x; k, l) :=& \sum_{\substack {n\leq x \\ n\equiv l \mod k}} \frac {\Lambda(n)}{ \log n}
\end{align*}
\end{dfn}
We have the corresponding analogues for $\delta_\pi(x; k, l_1, l_2)$, where the subscript is replaced by a different prime-counting function:
\begin{align*}
\delta_\psi(x; k, l_1, l_2):=&\psi(x; k, l_1)-\psi(x; k, l_2)\\
\delta_\vartheta(x; k, l_1, l_2):=&\vartheta(x; k, l_1) -\vartheta(x; k, l_2)  \\
\delta_\Pi(x; k, l_1, l_2):=&\Pi(x; k, l_1) -\Pi(x; k, l_2) 
\end{align*}
Further we define $w_f(T; k, l_1, l_2)$ to be the number of sign changes of $\delta_f(x; k, l_1, l_2)$ in the interval $[0, T]$ with fixed $k$, where $f \in \{\pi, \psi, \Pi, \theta\}$.

Since Chebyshev's original paper dealt with the case where each term in the sum contains a factor of $e^{-10x}$, we would able to form the {\em{mutatis mutadis}} definitions if we were to multiply a $e^{-nr}$ term to each term in the above sums:
\begin{align*}
\psi(x; k, l) \; \; \; \; \; {\rm to} \; \; \; &\sum_{n\equiv l \mod k} \Lambda(n)e^{-nr} \\
\Pi(x; k, l) \; \; \; \; \; {\rm to}  \; \; \; &\sum_{n\equiv l \mod k}\frac{ \Lambda(n)}{\log n}e^{-nr} \\
\vartheta(x; k, l) \; \; \; \; \; {\rm to} \; \; \; &\sum_{n\equiv l \mod k} \log{p}\, e^{-nr} \\
\pi(x; k, l) \; \; \; \; \; {\rm to} \; \; \;  &\sum_{n\equiv l \mod k} e^{-nr}\\
{\Li}(x) \; \; \; \; {\rm to} \; \; \;  &\int_0 ^\infty \frac{e^{-yr}}{\log y}\, dy 
\end{align*}

\begin{dfn}
The difference functions $\delta_f$ are replaced by $\Delta_F$'s:
\begin{align*}
\Delta_\psi(r; k, l_1, l_2):=&\sum_{n\equiv l_1 \mod k} \Lambda(n)e^{-nr} - \sum_{n\equiv l_2 \mod k} \Lambda(n)e^{-nr}\\
\Delta_\Pi(r; k, l_1, l_2):=&\sum_{n\equiv l_1 \mod k}\frac{ \Lambda(n)}{\log n}e^{-nr} - \sum_{n\equiv l_2 \mod k}\frac{ \Lambda(n)}{\log n}e^{-nr} \\
\Delta_\vartheta(r; k, l_1, l_2):=&\sum_{n\equiv l_1 \mod k} \log{p} e^{-nr}- \sum_{n\equiv l_2 \mod k} \log{p} e^{-nr} \\
\Delta_\pi(r; k, l_1, l_2):=& \sum_{n\equiv l_1 \mod k} e^{-nr} - \sum_{n\equiv l_2 \mod k} e^{-nr}  .
\end{align*}
Similarly  $w_f(T; k, l_1, l_2)$ is replaced by $W_F(T; k, l_1, l_2)$, for $F \in \{\psi, \Pi, \vartheta, \pi\}$.
\end{dfn}

\section{The Classical Tools}
The classical methods used to investigate the oscillatory properties of the functions $\delta_f(x; k, l_1, l_2)$ and $\Delta_F(x; k, l_1, l_2)$ for $F \in \{\psi, \Pi, \vartheta, \pi\}$ are inspired by the ones used to study the oscillatory term of the prime number theorem, namely $\pi(x) - \Li(x)$.
The primary tools are called the ``explicit formulas'', linking the functions $\pi(x; k, l)$ to the distribution of zeros of the Dirichlet-$L$ functions $L(s, \chi)$ in the critical strip, $0< \Re(s)<1$, for characters $\chi$ modulo $k$.

The asymptotic formula for $\pi(x; k, l)$ gives
\begin{equation}\label{asym}
\pi(x; k, 1)=\Pi(x; k, l)-\frac{N_k(l)}{\varphi(k)}\frac{x^{1/2}}{\log x}+o\Bigg(\frac{x^{1/2}}{\log x} \Bigg) \;\;\;\;\;(x \to \infty)
\end{equation}
where $N_k(l)$ denotes the number of incongruent solutions of the congruence $a^2 \equiv l \mod k$.

\begin{dfn}
Denote by $D_k$ and $C_k$ the sets of all characters and all non-principal characters modulo $k$, respectively. For $\chi \in D_k$, define
\[
\Psi(x; \chi):=\sum_{n\leq x}\Lambda(n)\chi(n).
\] 
\end{dfn}

Then it follows that
\begin{align*}
\varphi(k)\Pi(x; k, l) &{}= \varphi(k)\bigg(\frac{\psi(x; k, l)}{\log x}+\int_2^x\frac{\psi(t; k, l)}{t \log^2 t}\,dt \bigg)\\
&= \sum_{\chi \in D_k}\bar{\chi}(l)\bigg(\frac{\Psi(x; \chi)}{\log x}+\int_2^x\frac{\Psi(t; \chi)}{t \log^2 t}\,dt \bigg)
\end{align*}
and by equation~\eqref{asym} we have:
\begin{align*}
\varphi(k)\delta_\pi(x; k, l_1, l_2)&=\varphi(k)\bigg(\frac{\delta_\psi(x; k, l_1, l_2)}{\log x}+\int_2^x\frac{\delta_\psi(t; k, l_1, l_2)}{t \log^2 t}\, dt \bigg)\\
&\qquad-\big(N_k(l_1)-N_k(l_2) \big)\frac{x^{1/2}}{\log x}+o\bigg(\frac{x^{1/2}}{\log x}\bigg)\\
&=\sum_{\chi \in C_k}\big(\bar{\chi}(l_1)-\bar{\chi}(l_2)\big)\bigg(\frac{\Psi(x; \chi)}{\log x}+\int_2^x\frac{\Psi(t; \chi)}{t \log^2 t}\,dt\bigg)\\
&\qquad-\big(N_k(l_1)-N_k(l_2) \big)\frac{x^{1/2}}{\log x}+o\bigg(\frac{x^{1/2}}{\log x}\bigg) \;\;\;\; (x \to \infty).
\end{align*}
Furthermore, for any $\chi \in C_k$, the well-known explicit formula tells us that
\begin{equation}\label{Psi}
\Psi(x, \chi)=-\sum_{|\Im(\varrho)|\leq x}\frac{x^\varrho}{\varrho}+O(\log^2 x) \;\;\;\; (x \geq 2),
\end{equation}
where the sum runs over zeros $\varrho$ of $L(s, \chi)$ in the critical strip.
Now we see that the zeros of $L(s, \chi)$ play an important role in determining the distribution of primes to different moduli, and the zeros with the largest real part dominate the sum in equation \eqref{Psi}. Now we introduce a few conjectures that a handful subsequent results will require:

\begin{cnj}
[Generalized Riemann Hypothesis (GRH)]\label{GRH}
For any Dirichlet character $\chi$, all zeros of $L(s, \chi)$ inside the critical strip lie on the critical line $\sigma:=\Re(s)=\h$.
\end{cnj}
\noindent This of course is a generalization of the famous ``Riemann Hypothesis'', where we take $\chi$ to be the trivial character:
\begin{cnj}
[Riemann Hypothesis (RH)]\label{RH}
Inside the critical strip, the only zeros of the Riemann zeta function 
\[
\zeta(s):=\sum_{n=1}^\infty \frac{1}{n^s}
\]
satisfy $\sigma>0$ at $\sigma=\h$.
\end{cnj}

A basic tool for proving oscillation theorems is inspired by Landau's work~\cite{1905.Landau} on the location of singularities of the Mellin transforms of a non-negative function.  Suppose $f(x)$ is real-valued and non-negative for $x$ sufficiently large. Suppose also for some real numbers $\beta<\sigma$ that the Mellin transform
\[
g(s):=\int_1^\infty f(x)x^{-s-1}\,dx
\]
is analytic for $\Re(s)> \sigma$ and can be analytically continued to the real segment $(\beta, \sigma]$. Then $g(s)$ represents an analytic function in the half-plane $\Re(s)>\beta$.
\begin{xmp}
If $f(x)=\varphi(k)\delta_\psi(x; k, l_1, l_2)$ then
\[
g(s)=g(s; k, l_1, l_2)=-\frac{1}{s}\sum_{\chi \in C_k}\big(\bar{\chi}(l_1)-\bar{\chi}(l_2) \big)\frac{L'(s, \chi)}{L(s, \chi)}
\]
for $\Re(s)>1$, with the R.H.S. providing a meromorphic continuation of $g(s)$ to the whole complex plane.
\end{xmp}
\begin{rmk} Note that the poles of $g(s)$ above (except at $s=0$) are a subset of the zeros of the functions $L(s, \chi)$.  Also, $g(s)$ always has an infinite number of poles in the critical strip. Now assuming $g(s)$ with no real poles $s> \h$ and a pole $s_0$ with $\Re(s_0)>\h$, we take $\alpha$ satisfying $\h<\alpha<\Re(s_0)$ and put $f(x)=(-1)^n\delta_\psi(x; k, l_1, l_2)+c_1x^\alpha$ for some constant\footnote{Through out this paper, $c_i$ with $i\in \N$ shall always denote a calculable positive constant.} $c_1$ and $n \in \{0, 1 \}$.  The above discussion on Mellin transforms with different $n$ and $c_1$ yield that
\[
\limsup_{x \to \infty}\frac{\delta_\psi(x; k, l_1, l_2)}{x^\alpha}=+\infty, \;\;\;\;\;\;\; \liminf_{x \to \infty}\frac{\delta_\psi(x; k, l_1, l_2)}{x^\alpha}=-\infty 
\] 
\end{rmk}\section{Classical Results by Knapowski and Tur\'an, Serie I}\label{sec.Results1}
As mentioned in Section~\ref{intro section}, Knapowski and Tur\'an exhibited a  keen interest on this topic: they listed most of the problems in Section~\ref{intro section} and attempted to answer a few of them in their series of 15 papers.
Their investigation begins with the comparison of the progressions
\begin{equation*}
n\equiv 1 \mod k \;\;{\rm and } n\equiv l \mod k,\;\;\;\; {\rm where } 
\;\; l\not\equiv 1 \mod k.
\end{equation*}
First with 
\begin{equation}\label{II.2}
k=3, 4, 5, 6, 7, 8, 9, 10, 11, 12, 19, 24,
\end{equation}
which are in fact the first few numbers known to satisfy:
\begin{cnj}[Haselgrove Condition (HC) for the modulus $k$]\label{HC}
There is a function $Z(k)$ with $0< Z(k)\leq 1$ such that no $L(s, \chi)$ with $\chi \mod k$ vanishes for $0< \sigma < 1$, $|t|\leq Z(k)$, where $s=\sigma+it$, as usual.
\end{cnj}

\begin{dfn}
Define the iterated exponential and logarithmic functions by:
\begin{align*}
e_1(x):&=e^x, \;\;\;\;\;\;\;\;\, e_\nu(x):=e_{\nu-1}\big(\exp(x)\big)\\
\log_1(x):&=\log x, \;\;\;\;\; \log_\nu(x):=\log_{\nu-1}\big(\log(x)\big)
\end{align*}
\end{dfn}
\begin{thm}[\cite{1962.Knaposwki_2} Theorem 1.1] For any $k$ in \eqref{II.2}
\begin{align*}
\max_{T^{{1/3}}\leq x \leq T}\delta_\psi(x; k, 1, l)>& \sqrt{T}\exp\bigg{(} -41\frac{\log(T)\log_3{(T)}}{\log_2{(T)}}\bigg)\\
\min_{T^{{1/3}}\leq x \leq T}\delta_\psi(x; k, 1, l)<&- \sqrt{T}\exp\bigg{(} -41\frac{\log(T)\log_3{(T)}}{\log_2{(T)}}\bigg)
\end{align*}
\end{thm}
This essentially solves Problem~\ref{P:L-S} for $\delta_\psi$ in Section~\ref{intro section} with the $k$'s in equation \eqref{II.2}, in the case of $l_2=1$ at least.  Since C. L. Siegel proved~\cite{1945.Siegel} that for all $L(s, \chi)$ functions with primitive characters mod $k$ there is at least one zero $\varrho^*=\varrho^*(\chi)$ in the domain
\begin{equation}\label{Siegel}
\sigma \geq \h,\;\;\;\;\; |t|\leq \frac{c_2}{\log_3\big(k+e_3(1) \big)}, \;\;\;\;\;(s=\sigma +it)
\end{equation}
the above theorem follows at once from:
\begin{thm}[\cite{1962.Knaposwki_2} Theorem 1.2]\label{1.2 in II}
For a $k$ in \eqref{II.2} and a $\varrho_0=\beta_0+i\gamma_0$ with 
\begin{equation}\label{II 1.4}
\beta_0\geq \h,\;\;\;\; \gamma_0 > 0,
\end{equation}
where $\varrho_0$ is a zero of an $L(s, \chi^*)$ belonging to modulo $k$ with $\chi^*(l)\neq 1$ and $T>\max\big(c_3, e_2(10|\varrho_0|)\big)$, then the inequalities 
\begin{align*}
\max_{T^{{1/3}}\leq x \leq T}\delta_\psi(x; k, 1, l)>& T^{\beta_0}\exp\bigg{(} -41\frac{\log(T)\log_3{(T)}}{\log_2{(T)}}\bigg)\\
\min_{T^{1/3}\leq x \leq T}\delta_\psi(x; k, 1, l)< &-T^{\beta_0}\exp\bigg{(} -41\frac{\log(T)\log_3{(T)}}{\log_2{(T)}}\bigg)
\end{align*}
hold.
\end{thm}
The authors juxtapose a similar result:
\begin{thm}[\cite{1962.Knaposwki_2} Theorem 2.1]
For a $k$ in \eqref{II.2} and $T>c_4$ we have:
\begin{align*}
\max_{T^{1/3}\leq x \leq T}\delta_\Pi(x; k, 1, l)>& \sqrt{T}\exp\bigg{(} -41\frac{\log(T)\log_3{(T)}}{\log_2{(T)}}\bigg)\\
\min_{T^{1/3}\leq x \leq T}\delta_\Pi(x; k, 1, l)< &- \sqrt{T}\exp\bigg{(} -41\frac{\log(T)\log_3{(T)}}{\log_2{(T)}}\bigg)
\end{align*}
\end{thm}
and the above theorem essentially solves the case $l_2=1$ for the $k$'s in~\eqref{II.2}, for Problem~\ref{P:L-S} with $\delta_\Pi$. Now by~\eqref{Siegel} this theorem is an immediate consequence of: 
\begin{thm}[\cite{1962.Knaposwki_2} Theorem 2.2]\label{2.2 in II}
Let $k$ be a number in~\eqref{II.2}. If $\varrho_0$ with~\eqref{II 1.4} holds is a zero of an $L(s, \chi^*)$ belonging to modulo $k$, $\chi^*(l)\neq 1$ and $T>\max\big(c_5, e_2(10|\varrho_0|)\big)$, then the inequalities 
\begin{align*}
\max_{T^{1/3}\leq x \leq T}\delta_\Pi(x; k, 1, l)>& T^{\beta_0}\exp\bigg{(} -41\frac{\log(T)\log_3{(T)}}{\log_2{(T)}}\bigg)\\
\min_{T^{1/3}\leq x \leq T}\delta_\Pi(x; k, 1, l)<&-T^{\beta_0}\exp\bigg{(} -41\frac{\log(T)\log_3{(T)}}{\log_2{(T)}}\bigg)
\end{align*}
hold.
\end{thm}
Combining Theorems~\ref{1.2 in II} and~\ref{2.2 in II} (Theorems 1.2 and 2.2 in~\cite{1962.Knaposwki_2}) yields:
\begin{thm}[\cite{1962.Knaposwki_2} Theorem 3.1]\label{3.1 in II}
For a modulus $k$ satisfying the HC (Conjecture~\ref{HC}) and for a $\varrho_0$ with \eqref{II 1.4} holds, when 
\begin{equation*}
T>\max\bigg(c_6, e_2(10|\varrho_0|), e_2(k), e_2\Big(\frac{1}{Z(k)^3}\Big)\bigg)
\end{equation*}
we have the inequalities
\begin{align*}
\max_{T^{1/3}\leq x \leq T}\delta_\psi(x; k, 1, l)>& T^{\beta_0}\exp\bigg{(} -41\frac{\log(T)\log_3{(T)}}{\log_2{(T)}}\bigg)\\
\min_{T^{1/3}\leq x \leq T}\delta_\psi(x; k, 1, l)<&-T^{\beta_0}\exp\bigg{(} -41\frac{\log(T)\log_3{(T)}}{\log_2{(T)}}\bigg)
\end{align*}
and further
\begin{align*}
\max_{T^{1/3}\leq x \leq T}\delta_\pi(x; k, 1, l)>& T^{\beta_0}\exp\bigg{(} -41\frac{\log(T)\log_3{(T)}}{\log_2{(T)}}\bigg)\\
\min_{T^{1/3}\leq x \leq T}\delta_\pi(x; k, 1, l)<&-T^{\beta_0}\exp\bigg{(} -41\frac{\log(T)\log_3{(T)}}{\log_2{(T)}}\bigg)
\end{align*}
\end{thm}
hold.\\
This concludes their answer to Problems~\ref{P:L-S} for $\delta_\pi, \delta_\Pi$ and $\delta_\psi$, at least for the case $l_2=1$.
\vskip12pt
The authors move on to other problems and their variations, where they first give, as a consequence of Theorem~\ref{3.1 in II} by taking $\varrho^*$ satisfying the condition of Siegel's Theorem~\eqref{Siegel}:
\begin{thm}[\cite{1962.Knaposwki_2} Theorem 4.1]
In the interval 
\[
0<x<\max \bigg(c_7, e_2(k), e_2\Big(\frac{1}{Z(k)^3} \Big) \bigg)
\]
the functions $\delta_\psi(x; k, 1, l)$ and $\delta_\Pi(x; k, 1, l)$ certainly change their sign, when $k$ satisfies the HC (Conjecture~\ref{HC}). 
Here
\begin{equation*}
c_7=\max\big(c_6, e_2(10(1+c_2))\big)
\end{equation*}
\end{thm}
The above theorem gives answers to Problem~\ref{P:F-S} for $\delta_\psi$ and $\delta_\Pi$, and the authors conjecture the ``best'' interval is 
\[
0<x<\exp(c_8k)
\]
Then they appeal to some answers for Problem~\ref{P:A-B} regarding $w_\psi$ and $w_\Pi$, giving:
\begin{thm}[\cite{1962.Knaposwki_2} Theorem 4.2]
If for a $k$ holing the HC (Conjecture~\ref{HC}) and with 
\[
T>\exp\Bigg( c_9 \bigg(\exp(k)+\exp\Big(\frac{1}{Z(k)^3} \Big) \bigg)^2\Bigg)
\]
the inequalities
\begin{align*}
w_\psi(T; k, 1, l_1) &> \frac{1}{8\log 3}\log_2{T} \\
w_\Pi(T; k, 1, l_1) &> \frac{1}{8\log 3}\log_2{T}
\end{align*}
hold.
\end{thm}
As an obvious consequence of Theorem~\ref{3.1 in II}, they assert,
\begin{thm}[\cite{1962.Knaposwki_2} Theorem 4.3]
Let $L(s, \chi^*)$ be an arbitrary $L$-function mod $k$ ($k$ holding HC (Conjecture~\ref{HC})), and for
\begin{equation*}
T> \max\bigg(c_6, e_2(k), e_2\Big(\frac{1}{Z(k)^3}\Big)\bigg),
\end{equation*}
if $l$ is such that $\chi^*(l)\neq 1$ , then $L(s, \chi)$ does not vanish in the domain 
$$
\sigma \geq 41\frac{\log_3{T}}{\log_2{T}} +\frac{1}{\log T}\max_{T^{{1/3}} \leq x \leq T}\log \delta_\psi(x; k, 1, l)
$$
$$
|t|\leq \frac{1}{10} \log_2{T} -1.
$$
\end{thm}
\vskip12pt
Then they give partial answers to Problem~\ref{P:L-G}:
\begin{thm}[\cite{1962.Knaposwki_2} Theorem 5.1]
If $k$ is one of the moduli \eqref{II.2} then for $T> c_{10}$, we have the inequalities
\begin{align}
\label{II.5.1.a}
\max_{\exp(\log_3^{1/130}T ) \leq x \leq T} \frac{\delta_\pi(x; k, 1, l)}{\bigg( \frac{\sqrt{x}}{\log x}\bigg)}  &> \frac{1}{100} \log_5 T\\
\label{II.5.1.b}
\min_{\exp(\log_3^{1/130}T ) \leq x \leq T} \frac{\delta_\pi(x; k, 1, l)}{\bigg( \frac{\sqrt{x}}{\log x}\bigg)}  &< - \frac{1}{100} \log_5 T
\end{align}
\end{thm}
\begin{thm}[\cite{1962.Knaposwki_2} Theorem 5.2]
If HC (Conjecture~\ref{HC}) holds for a $k$ and
\begin{equation} 
\label{II.5.2}
T>\max\bigg( e_5(c_{11} k), e_2\Big(\frac{1}{Z(k)^3}\Big)\bigg)
\end{equation}
then the inequalities~\ref{II.5.1.a} and~\ref{II.5.1.b} hold.
\end{thm}
and further for Problem~\ref{P:F-S}, there is
\begin{thm}[\cite{1962.Knaposwki_2} Theorem 5.3]
If HC (Conjecture~\ref{HC}) holds for a $k$ then the interval 
\begin{equation*} 
1\leq x \leq \max \bigg( e_5(c_{11} k), e_2\Big(\frac{1}{Z(k)^3}\Big)\bigg)
\end{equation*}
contains at least a zero of $\delta_\pi(x; k, 1, l)$.
\end{thm}

\vskip12pt
As for Problem~\ref{P:P-R}, they gave:
\begin{thm}[\cite{1962.Knaposwki_2} Theorem 5.4]
If HC (Conjecture~\ref{HC}) holds for a $k$ and for $T$ with \eqref{II.5.2}, the inequalities
\begin{align*}
\max_{\exp(\log_3^{1/130}T )} \frac{\log x}{\sqrt{x}} \bigg\{ \pi(x; k, 1) - \frac{1}{\varphi(x)}\pi(x)
\bigg\} &> \frac{1}{200} \log_5{T}\\
\min_{\exp(\log_3^{1/130}T )} \frac{\log x}{\sqrt{x}} \bigg\{ \pi(x; k, 1) - \frac{1}{\varphi(x)}\pi(x)
\bigg\} &<- \frac{1}{200} \log_5{T}.
\end{align*}
hold.
\end{thm}
\vskip12pt
Revisiting Problem~\ref{P:A-B}, the authors present:
\begin{thm}[\cite{1962.Knaposwki_3} Theorem 1.1]\label{1.1 in III}
For $T>c_{12}$ we have for the moduli $k$ in \eqref{II.2} the inequality 
\begin{equation*}
w_\pi(T; k, 1, l)> c_{13} \log_4 T
\end{equation*}
holds.
\end{thm}
and more generally,
\begin{thm}[\cite{1962.Knaposwki_3} Theorem 1.2]\label{1.2 in III}
If $k$ satisfies the HC (Conjecture~\ref{HC}) holds then for 
\begin{equation*}
T>\max \bigg(e_4(k^{c_{14}}), e_2\Big(\frac{2}{Z(k)^3} \Big)\bigg)
\end{equation*}
we have the inequality 
\begin{equation*}
w_\pi(T; k, 1, l)>k^{-c_{14}}\log_4 T.
\end{equation*}
\end{thm}
As a consequence they show
\begin{thm}[\cite{1962.Knaposwki_3} Theorem 1.3]
If for a $k$ holding the HC (Conjecture~\ref{HC})  then in the interval
\[
0 < x \leq \max  \bigg(e_4(k^{c_{14}}), e_2\Big(\frac{2}{Z(k)^3} \Big)\bigg)
\]
there exists at least one $x$ such that 
$ \delta_\pi (T; k, 1, l)=0 $ changes its sign
for all $l\not\equiv 1 \mod k$.
\end{thm}

They also prove the analogous theorems of Theorem~\ref{1.1 in III} and Theorem~\ref{1.2 in III} with the following definition, contributing to Problem~\ref{P:P-R}:
\begin{dfn}
For 
\begin{equation}\label{difference}
\pi(x; k, 1)- \frac{1}{\varphi(k)-1}\sum_{\substack{(l, k) =1 \\ l \neq 1 \\ l}} \pi(x; k, l)
\end{equation}
we denote the number of sign-changes in this function for $x \in (0,  T] $ by $S_k(T)$
\end{dfn}
\begin{thm}[\cite{1962.Knaposwki_3} Theorem 1.4]
If $k$ satisfies the HC (Conjecture~\ref{HC}) then for 
$$
T>\max \bigg(e_4(k^{c_{14}}), e_2\Big(\frac{2}{Z(k)^3} \Big)\bigg)
$$
the inequality 
$$
S_k(T)> k^{-c_{14}}\log_4{T}
$$
holds, and the same result holds if we changed formula~\eqref{difference} to 
\[\pi(x; k, 1)-\frac{1}{\varphi(x)}\pi(x)
\;\;\;\;{\rm and}\;\;\;\;\pi(x; k, 1)-\frac{1}{\varphi(x)}\Li(x)
\]
\end{thm}
\vskip12pt
The authors again revisit Problems~\ref{P:A-B} and~\ref{P:F-S}, armed with the above theorems, and assuming that equations \eqref{l_1} and \eqref{l_2} having exactly the same number of solutions, whose significance we speculated in Remark~\ref{l_1 and l_2}:
\begin{thm}[\cite{1962.Knaposwki_3} Theorem 2.1]\label{2.1 of III}
For the $k$'s in \eqref{II.2} and $l$'s satisfying the condition \eqref{l_1} and \eqref{l_2}  having the same number of solutions with $l_2=1$ then for $T> c_{14}$ the inequalities
\begin{align*}
\max_{T^{1/3}\leq x \leq T}\delta_\pi(x; k, 1, l)>& \sqrt{T}\exp\bigg{(} -41\frac{\log(T)\log_3{(T)}}{\log_2{(T)}}\bigg)\\
\min_{T^{1/3}\leq x \leq T}\delta_\pi(x; k, 1, l)<&-\sqrt{T}\exp\bigg{(} -41\frac{\log(T)\log_3{(T)}}{\log_2{(T)}}\bigg)
\end{align*}
hold
\end{thm}

which is a special case of:
\begin{thm}[\cite{1962.Knaposwki_3} Theorem 2.2]\label{2.2 of III}
For the modulus $k$ in \eqref{II.2}, for a $l$ satisfying~\eqref{l_1} and~\eqref{l_2} with $l_2=1$, of $\varrho_0=\beta_0+i\gamma_0$ with $\beta_0 \geq \h$ such that $L( \varrho_0, \chi) = 0$ with $\chi(l)\neq 1$, then we have for
\[
T> \max\big(c_{15}, e_2(10|\varrho_0|)\big)
\] 
the inequalities 
\begin{align*}
\max_{T^{1/3}\leq x \leq T}\delta_\pi(x; k, 1, l)>& T^{\beta_0}\exp\bigg{(} -41\frac{\log(T)\log_3{(T)}}{\log_2{(T)}}\bigg)\\
\min_{T^{1/3}\leq x \leq T}\delta_\pi(x; k, 1, l)<&-T^{\beta_0}\exp\bigg{(} -41\frac{\log(T)\log_3{(T)}}{\log_2{(T)}}\bigg)
\end{align*}
hold.
\end{thm}
due to Siegel's Theorem~\eqref{Siegel}, as before.
\begin{thm}[\cite{1962.Knaposwki_3} Theorem 3.1]
If for a $k$ the HC (Conjecture~\ref{HC}) holds and $l$ satisfies \eqref{l_1} and \eqref{l_2} with $l_2=1$ then for
\[
T>\max\bigg(c_{16}, e_2(k), e_2\Big(\frac{1}{Z(k)^3} \Big)\bigg),
\]
the inequalities in Theorem~\ref{2.1 of III} hold.
\end{thm}

\begin{thm}[\cite{1962.Knaposwki_3} Theorem 3.2]
If for a $k$ the HC (Conjecture~\ref{HC}) holds and $l$ satisfting~\eqref{l_1} and~\eqref{l_2} with $l_2=1$, and if further $\varrho = \beta_0 + i\gamma$, $\beta_0 \geq \h$ is a zero for an $L(s, \chi)$ with $\chi(l) \neq 1$, then for 
\begin{equation*}
T>\max \bigg( c_{16}, e_2(k), e_2\Big(\frac{1}{Z(k)^3} \Big), e_2(10|\varrho|) \bigg),
\end{equation*}
the inequalities in Theorem~\ref{2.2 of III} hold.
\end{thm}

Now turning in to Problem~\ref{P:A-B} again:
\begin{thm}[\cite{1962.Knaposwki_3} Theorem 3.3]
For $T>c_1$ and $k$'s in the moduli \eqref{II.2} and $l's$ satisfying \eqref{l_1} and \eqref{l_2}, then the inequality
\[
w_\pi(T; k, 1, l)> c_{17} \log_2 T
\]
holds.
\end{thm}
and 
\begin{thm}[\cite{1962.Knaposwki_3} Theorem 3.4]
If for a $k$ the HC (Conjecture~\ref{HC}) holds and 
\[
T>\max\bigg(c_{18}, e_2(2k), e_2\Big(\frac{2}{Z(k)^3}\Big) \bigg)
\]
$l$ satisfies \eqref{l_1} and \eqref{l_2} then 
\[
w_\pi (T; K, 1, l) > c_{17} \log_2 T
\]
\end{thm}

The authors then delve into the general case of Problems~\ref{P:L-B} and~\ref{P:A-B} with $k = 8$ and $5$.
\begin{thm}[\cite{1963.Knaposwki_4} Theorem 1.1]
For $T>c_{18}$ and for all pairs $l_1$ and $l_2 $ with $l_1\neq l_2$ among the numbers 3, 5, 7 $\mod 8$, we have
\begin{equation*}
\max_{T^{{1/3}}\leq x \leq T}\delta_\pi(x; 8, l_1, l_2) > \sqrt{T} \Big(-23 \frac{\log T \log_3 T}{\log_2 T}\Big)
\end{equation*}
\end{thm}

\begin{thm}[\cite{1963.Knaposwki_4} Theorem 1.2]\label{1.2 of IV}
For $T>c_{18}$, the inequality 
\begin{equation*}
w_\pi(T; 8, l_1, l_2)>c_{19}\log_2 T
\end{equation*}
holds if only $l_1 \neq l_2$ among $3, 5, 7$.
\end{thm}

Since the congruence
\[
x^2\equiv l \mod 8, \; \; \; \; \; l\not\equiv 1 \mod 8
\]
is not solvable, it implies that Theorem~\ref{1.2 of IV} is a consequence of 

\begin{thm}[\cite{1963.Knaposwki_4} Theorem 2.1]
For $T>c_{20}$ and all pairs $l_1 \neq l_2$ among the numbers $3, 5, 7$ we have
\[
\max_{T^{{1/3}}\leq x \leq T}\delta_\Pi(x; 8, l_1, l_2) >\sqrt{T} \exp\Big(-23 \frac{\log T \log_3 T}{\log_2 T}\Big)
\]
\end{thm}
and with slight modifications we obtain:
\begin{thm}[\cite{1963.Knaposwki_4} Theorem 2.2]
For $T>c_{20}$ and all pairs $l_1 \neq l_2$ among the numbers $3, 5, 7$ we have
\[
\max_{T^{{1/3}}\leq x \lq T}\delta_\psi(x; 8, l_1, l_2) >\sqrt{T} \exp\Big(-23 \frac{\log T \log_3 T}{\log_2 T} \Big)
\]
\end{thm}
As an corollary we have:
\begin{thm}[\cite{1963.Knaposwki_4} Theorem 2.3]
For $T>c_{20}$ and all pairs $l_1 \neq l_2$ among the numbers $3, 5, 7$ we have
\begin{align*}
w_\psi(T; 8, l_1, l_2)> &\log_2 T \\
w_\Pi(T; 8, l_1, l_2)> &\log_2 T
\end{align*}
\end{thm}

continuing the study of the general cases, this time assuming  ``finite'' GRH (Conjecture~\ref{GRH}) Conjecture:
Problem~\ref{P:L-S} for $\delta_\psi(x; k, l_1, l_2)$ and $\delta_\Pi(x; k, l_1, l_2)$:
\begin{thm}[\cite{1963.Knaposwki_5} Theorem 1.1]\label{1.1 of V}
Supposing the truth of the ``finite'' GRH (Conjecture~\ref{GRH}), which says no $L(s, \chi)$ vanishes for a sufficiently large $c_{21} \geq 1$  
\begin{equation}\label{V 1.3}
\sigma >\h,  \; \; |t| \leq c_{21}k^{10},
\end{equation}
moreover also for 
\begin{equation}\label{V 1.4}
\sigma = \h, \;\; |t|\leq A(k)
\end{equation}
with $A(k)$ positive, $c_{22}$ sufficiently large, for
\begin{equation}\label{V 1.5}
T>\max\Bigg\{e_2(c_{22}k^{20}), \exp\bigg( 2\exp\Big(\frac{1}{A(k)^3}\Big)+c_{22}k^{20}\bigg)\Bigg\}
\end{equation}
we have for $l_1\neq l_2$ the inequalities:
\begin{align*}
\max_{T^{{1/3}}\leq x \leq T}\delta_\psi(x; k, l_1, l_2)>& \sqrt{T}\exp\bigg(-44\frac{\log T \log_3 T}{\log_2 T}\bigg) \\
\max_{T^{{1/3}}\leq x \leq T}\delta_\Pi(x; k, l_1, l_2)>& \sqrt{T}\exp\bigg(-44\frac{\log T \log_3 T}{\log_2 T}\bigg)
\end{align*}
\end{thm}

\begin{thm}[\cite{1963.Knaposwki_5} Theorem 1.2]
By the above theorem, both of $\delta_\psi(x; k, l_1, l_2)$ and  $\delta_\Pi(x; k, l_1, l_2)$ have a sign change in the interval $[T^{{1/3}},  T]$ whenever $T$ satisfies \eqref{V 1.5}, then we get at once:\\
For 
\begin{equation*}
T>\max\Bigg\{\exp\big(9\exp(2c_{22} k^{20}\big), \exp\bigg( 72\exp\Big(\frac{2}{A(k)^3}\Big)+18c_{22}^2k^{40}\bigg)\Bigg\}
\end{equation*} 
the inequalities 
\begin{align*}
w_\psi(T; k, l_1, l_2)&> \frac{\log_2 T}{2\log 3} \\
w_\Pi(T; k, l_1, l_2) &> \frac{\log_2 T}{2\log 3}
\end{align*}
hold.
\end{thm}

\begin{rmk}
We note that if $l_1$ and $l_2$ are such that none of \eqref{l_1} and \eqref{l_2} are solvable, then it follows from Theorem~\ref{1.1 of V} (with $c_{22}$ being replaced by a larger constant), that for
\[
T>\max\Bigg\{e_2\big(c_{23}k^{20}\big), \exp\bigg(2\exp\bigg(\frac{1}{A(k)^3} \bigg)+c_{24}k^{40}\bigg)\Bigg\}
\] 
the inequality
\[
\max_{T^{{1/3}}\leq x \leq T}\delta_\pi(x; k, l_1, l_2)>\sqrt{T}\exp\bigg(-45\frac{\log{T}\log_3{T}}{\log_2{T}} \bigg)
\]
holds.
\end{rmk}

Returning to Problem~\ref{P:P-R} with slight variations in the question:
\begin{thm}[\cite{1963.Knaposwki_5} Theorem 3.1]
Supposing the truth of finite GRH (Conjecture~\ref{GRH}), we have for each $(l, k)=1$ and
\begin{equation*}
T>\max\Bigg\{e_2(c_{22}k^{20}), \exp\bigg(2\exp\Big(\frac{1}{A(k)^3}\Big)+c_{23}k^{20}\bigg)\Bigg\}
\end{equation*}
both the inequalities
\begin{align*}
\max_{T^{1/3}\leq x \leq T}\bigg\{\Pi(x, k, l)-\frac{1}{\varphi(k)}\Pi(x)\bigg\}>&\sqrt{T}\exp\Big(-44\frac{\log T\log_3 T}{\log_2 T}\Big) \\
\max_{T^{1/3}\leq x \leq T}\bigg\{\Pi(x, k, l)-\frac{1}{\varphi(k)}\Pi(x)\bigg\}<&-\sqrt{T}\exp\Big(-44\frac{\log T\log_3 T}{\log_2 T}\Big)
\end{align*}
hold, and the same hold if we replace $\Pi$ by $\psi$.
\end{thm}

\begin{rmk}
The analogous statements also hold for if we change
\[
\Pi(x, k, l)-\frac{1}{\varphi(k)}\Pi(x)
\]
to
\[
\Pi(x; k, l)-\frac{1}{\varphi(k)}\Li(x)  \;\;\;\;{\rm{and}}\;\;\;\; \psi(x; k. l)-\frac{1}{\varphi(k)}\Li(x)
\]
in the above theorem.  However main difficulties occur when trying to prove that for all $(l, k)=1$ the function
\[
\pi(x; k, l)-\frac{1}{\varphi(k)}\Li(x)
\]
changes sign infinitely often.  We speculate that this is plausible if only the congruence $x^2 \equiv l \mod k$ is not solvable, as we succeeded in proving similar results in the other cases.
\end{rmk}

\begin{thm}[\cite{1963.Knaposwki_6} Theorem 1.1]
If for a $k$ the assertions \eqref{V 1.3}  and \eqref{V 1.4} hold, then for 
\begin{equation}
T>\max\Bigg\{e_2(c_{24}k^{20}), \exp\bigg(2\exp\Big(\frac{1}{A(k)^3}\Big)+c_{24}k^{20}\bigg)\Bigg\}
\end{equation} 
and all $(l_1, l_2)$ pairs of two squares or two non-squares mod $k$, both the following inequalities hold:
\begin{align}
\max_{T^{{1/3}}\leq x \leq T}\delta_\pi(x; k, l_1, l_2)>&\sqrt{T}\exp\Big(-44\frac{\log T\log_3 T}{\log_2 T}\Big) \\
\max_{T^{{1/3}}\leq x \leq T}\delta_\pi(x; k, l_2, l_1)<-&\sqrt{T}\exp\Big(-44\frac{\log T\log_3 T}{\log_2 T}\Big) 
\end{align}
\end{thm}

now they examine how $\delta_\psi(x; k, l_1, l_2)$ changes its signs infinitely often:
\begin{thm}[\cite{1963.Knaposwki_7} Theorem 1.1]\label{VII 1.1}
 Answers Problem~\ref{P:S-C} for $\delta_\psi(x; k, l_1, l_2)$:
Under GRH (Conjecture~\ref{GRH}) for $\chi$ mod $k$, each function $\delta_\psi(x; k, l_1, l_2)$ with $l_1 \neq l_2$ changes its sign infinitely often for $1\leq x < +\infty$
\end{thm}

\begin{thm}[\cite{1963.Knaposwki_7} Theorem 1.2]\label{VII 1.2}
Regarding Problem~\ref{P:F-S} for the case of $\delta_\psi(x; k, l_1, l_2)$:
First sign change:
For all $k$'s satisfying the HC (Conjecture~\ref{HC}), all functions $\delta_\psi(x; k, l_1, l_2)$ change their sign in the interval 
\begin{equation*}
1\leq x \leq \max\bigg( e_2(k^{c_{25}}), e_2\Big( \frac{2}{Z(k)^3} \Big) \bigg)
\end{equation*}
with a sufficiently large $c_{25}$
\end{thm}

Both of Theorems~\ref{VII 1.1} and~\ref{VII 1.2} easily follows from
\begin{thm}[\cite{1963.Knaposwki_7} Theorem 1.3]
For the $k$'s satisfying the HC (Conjecture~\ref{HC}), all functions $\delta_\psi(x; k, l_1, l_2)$ change their sign in the interval 
\begin{equation*}
\omega \leq x \leq e^{2\sqrt{\omega}}
\end{equation*}
only if
\begin{equation*}
\omega \geq \max \bigg( \exp(k^{c_{26}}), \exp\Big(\frac{2}{Z(k)^3}\Big)\bigg)
\end{equation*}
for a sufficiently large $c_{26}$.
\end{thm}

Finally we have some unconditional results at the end of serie 1, for $k=8$
\begin{thm}[\cite{1963.Knaposwki_8} Theorem 1.1]\label{VIII 1.1}
If $0<\delta< c_{27}$, then for $l_1\not\equiv l_2 \not\equiv 1 \mod 8$ the inequality
\begin{equation}
\max_{\delta\leq x \leq \delta^{\frac 1 3}} 
\Delta_\vartheta(x; 8, l_1, l_2)> \frac{1} {\sqrt\delta} \exp \Bigg(-22\frac {\log\big( 1/\delta\big) \log_3\big( 1/\delta)} {\log_2\big( 1/\delta \big)} \Bigg)
\end{equation}
holds unconditionally, and since $l_1$ and $l_2$ can be interchanged,
\begin{equation*}
\max_{\delta\leq x \leq \delta^{\frac 1 3}} 
\Delta_\vartheta(x; 8, l_1, l_2)<- \frac{1} {\sqrt\delta} \exp \Bigg(-22\frac {\log\big( 1/\delta\big) \log_3\big( 1/\delta)} {\log_2\big( 1/\delta \big)} \Bigg)
\end{equation*}
also holds.
\end{thm}
For the case $l_1=1$ they show
\begin{thm}[\cite{1963.Knaposwki_8} Theorem 1.2]\label{VIII 1.2}
If for an $l \not\equiv 1 \mod 8$ 
\[
\lim_{x\to +0}\Delta_\vartheta(x; k, 1, l) = -\infty
\]
then no $L(s, \chi)$-function mod 8 with $\chi(x) \neq 1$ can vanish for $\sigma > \h$
\end{thm}
Further they prove:
\begin{thm}[\cite{1963.Knaposwki_8} Theorem 1.3]\label{VIII 1.3}
If no $L(s, \chi)$ functions $\mod 8$ with  $\chi \in C_k$ vanish for $\sigma > \h$, then for all $l\not\equiv 1\mod 8$ we have
If for an $l \not\equiv 1 \mod 8$ 
\[
\lim_{x\to  +0}\Delta_\vartheta(x; k, 1, l) = -\infty
\]
\end{thm}

\begin{rmk}
Theorem~\ref{VIII 1.2} can be shown by mimicking Landau's argument~\cite{1918.Landau.1} with slight modifications, and Theorem~\ref{VIII 1.3} by Hardy-Littlewood-Landau's argument~\cite{1918.Landau.1}~\cite{1918.Landau.2}, and \cite{1918.Littlewood}. On the other hand, of course Theorem~\ref{VIII 1.1} more difficult to show. Since the congruences
\[
x^2\equiv l \mod 8 \;\;\;\;l=3, 5, 7
\]  
and hence
\[
\max_{\delta \leq x \leq \delta^{\frac 1 3}}\sum_{\substack{p, \nu \\ \nu \geq 3}}\log p\cdot\exp(-p^\nu x)=O\bigg(\frac{1}{\delta^{{1/3}}}\log^2\frac{1}{\delta} \bigg)
\]
Theorem~\ref{VIII 1.1} is equivalent to the inequality
\begin{equation*}
\max_{\delta\leq x \leq \delta^{\frac 1 3}} 
\Delta_\vartheta(x; 8, l_1, l_2)> \frac{1} {\sqrt\delta} \exp \Bigg(-22\frac {\log\big( 1/\delta\big) \log_3\big( 1/\delta)} {\log_2\big( 1/\delta \big)} \Bigg)
\end{equation*}
\end{rmk}

\begin{thm}[\cite{1963.Knaposwki_8} Theorem 1.4]
If  $0<\delta< c_{28}\leq 1$, then for $l\not\equiv 1 \mod8 $ the following inequalities hold:
\begin{align*}
\max_{\delta\leq x \leq \delta^{\frac 1 3}}  \Delta_{\psi}(x; 8, 1, l)&>\frac{1}{\sqrt\delta}\exp\ \Bigg(-22\frac {\log\big( 1/\delta\big) \log_3\big( 1/\delta)} {\log_2\big( 1/\delta \big)} \Bigg) \\
\min_{\delta\leq x \leq \delta^{\frac 1 3}}  \Delta_{\psi}(x; 8, 1, l)&<-\frac{1}{\sqrt\delta}\exp\ \Bigg(-22\frac {\log\big( 1/\delta\big) \log_3\big( 1/\delta)} {\log_2\big( 1/\delta \big)} \Bigg)
\end{align*}
\end{thm}


\section{Classical Results by Knapowski and Tur\'an, Serie II}\label{sec.Results2}
\begin{thm}[\cite{1964.Turan_1}]
Let $k$ fulfill the HC (Conjecture~\ref{HC}), $(l, k) =1$, and let $\varrho=\beta+i\gamma$ be an zero of an $L(s, \chi)$ with $\chi(l)\neq 1$and $\beta\geq \h$. Then with a sufficiently large $c_{29}$ for 
\begin{equation*}
T>\max\bigg(c_{29}, e_2(k), \exp\Big(\frac{1}{Z(k)}\Big), e_2(|\varrho|)\bigg)
\end{equation*}
with suitable $U_1$, $U_2$, $U_3$, $U_4$ satisfying
\begin{align*}
T\exp\big(-\log^{11/12}T\big)\leq U_1 <U_2\leq T \\
T\exp\big(-\log^{11/12}T\big)\leq U_3 <U_4\leq T
\end{align*}
the inequalities 
\begin{align*}
\sum_{\substack{n\equiv 1 \mod k\\U_1\leq n \leq U_2}}\Lambda(n)-\sum_{\substack{n\equiv l \mod k\\U_1\leq n \leq U_2}}\Lambda(n) &\geq T^{\beta}\exp\big( -\log^{11/12}T\big)\\
\sum_{\substack{n\equiv 1 \mod k\\U_1\leq n \leq U_2}}\Lambda(n)-\sum_{\substack{n\equiv l \mod k\\U_3\leq n \leq U_4}}\Lambda(n) &\leq T^{\beta}\exp\big( -\log^{11/12}T\big)
\end{align*}
hold.
\end{thm}

\begin{dfn}
To study primes in different modulo, we adapt the following notation:
\[
\varepsilon(k; p, l_1, l_2):=\begin{cases} 1&\;\;\; {\rm if}\;\;\; p\equiv l_1 \mod k\\
                                          -1&\;\;\; {\rm if}\;\;\; p\equiv l_2 \mod k\\
                                           0&\;\;\; {\rm otherwise}
                                           \end{cases}
\]
\end{dfn}

\begin{xmp}
Hardy-Littlewood-Landau's argument~\cite{1918.Landau.1}~\cite{1918.Landau.2}, and~\cite{1918.Littlewood} gave (with abundant numerical data as well) that the relation:
\[
\lim_{x \to \infty}\sum_p\varepsilon(8; p, 1, l)\log(p)\exp\bigg(-\frac{p}{x}\bigg)=-\infty \;\;\;\;(l= 3, 5, 7)
\]
holds if and only if no $L(s, \chi)$ having $\mod 8$ with $\chi \in C_k$ vanishes for $\sigma > \h$ 
\end{xmp}
A few main results when $\exp\big(-\frac{p}{x}\big)$ is replaced by $\exp\Big( -\frac{1}{r(x)}\log^2\big(\frac p x \big)\Big)$ with suitable (``small'') $r(x)$:
\begin{thm}[\cite{1964.Turan_2} Theorem I]\label{IIb 1}
For any fixed $k$ satisfies the HC (Conjecture~\ref{HC}) and for all quadratic non-residues $l \mod k$, $(l, k) =1$, the relation
\begin{equation*}
\lim_{x \to \infty} \sum_p{\varepsilon(k; p, l, 1) \log{p}\cdot  \exp\bigg( -\frac{1}{r(x)}\log^2\Big(\frac p x \Big)\bigg)} = +\infty
\end{equation*}
for every $r(x)$ satisfying $0 < r(x) \leq \log x$ is valid if and only if none of the $L$-functions $\mod k$, with $\chi \in C_k$ vanishes for $\sigma > \h$
\end{thm}
which is a special case of:
\begin{thm}[\cite{1964.Turan_2} Theorem II]\label{IIb 2} 
For any fixed $k$ satisfies the HC (Conjecture~\ref{HC}) and for all quadratic non-residues $l \mod k$, $(l, k) =1$, the relation
\begin{equation*}
\lim_{x \to \infty} \sum_p{\varepsilon(k; p, l, 1) \log{p} \cdot  \exp\bigg( -\frac{1}{r(x)}\log^2\Big(\frac p x \Big)\bigg)} = +\infty
\end{equation*}
for every $r(x)$ satisfying $r_0< r(x) \leq \log x$ holds if and only if none of $L(s, \chi) \mod k$, with $\chi(l)\neq 1$ vanishes for $\sigma > \h$
\end{thm}
To deduce Theorem~\ref{IIb 1} from Theorem~\ref{IIb 2} we only have to note that for a character $\chi^*$, all non-residues $l$, $\chi^*(l)=1$, then $\chi^*$ is principle.

\begin{thm}[\cite{1964.Turan_2} Theorem III]\label{IIb 3}  
Assume $E(k) \leq {\sqrt{\log{k}}}/{k}$, if for a $k$ satisfying the HC (Conjecture~\ref{HC}) and a prescribed quadratic non-residue $l$, no $L(s, \chi)$ with $\chi(l) \neq 1$ vanishes for $\sigma > \h$, then for suitable $c_{30}, c_{31}, c_{32}$ and
\[
r_0 = c_{30} \frac{\log k}{E(k)^2}
\] 
the inequality
\begin{equation*}
\sum_p \varepsilon(k; p, l, 1) \log{p}  \exp \bigg( -\frac{1}{r(x)}\log^2\Big(\frac{p}{x} \Big)\bigg)> c_{31} \sqrt{x}
\end{equation*}
holds whenever $r_0< r \leq \log  x$ and $x > c_{32} k^{50}.$\\
As the contribution of primes $p$ with
\[
p>x\exp\big(10\sqrt{r\log x} \big)\;\;\;\;{\rm and}\;\;\;\;p<x\exp\big(-10\sqrt{r\log x} \big)
\]
is $o\big(\sqrt{x}\big)$, this theorem asserts under the given circumstances the preponderance of primes $\equiv l \mod k$ over those $\equiv 1 \mod k$ in the interval $\Big(x\exp(-10\sqrt{r\log x}), x\exp(10\sqrt{r\log x})\Big)$.
\end{thm}

\begin{thm}[\cite{1964.Turan_2} Theorem IV]  
Assume $E(k) \leq {\sqrt{\log{k}}}/{k}$, and if for a $k$ satisfying the HC (Conjecture~\ref{HC}) and a quadratic non-residue $l$ there exists an $L(s, \chi)$ with $\chi(k) \neq 1$ such that 
\begin{equation}\label{IIb 5.4}
L(\varrho_0, \chi) =0, \;\;\;\;\; \varrho_0=\beta + i\gamma, \;\;\;\;\; \beta>\h, \;\;\;\;\; \gamma>0
\end{equation} 
then for all $T$ with 
\begin{equation}\label{IIb 5.5}
T>\max \bigg(c_{33}, \exp\Big( \pi^7E(k)^{-7}\Big), \exp\Big( \exp(k)\Big), \exp\Big( \big( \frac{4+\gamma^2}{\beta -\h}\big) ^{21}\Big) \bigg)
\end{equation}
then there exist integers $r_1$ and $r_2$ with
\begin{equation*}
2\log^{5/7}{T} - 4\log^{4/7}{T}\leq r_1, r_2 \leq 2\log^{5/7}{T} - 4\log^{4/7}{T}
\end{equation*}
and $x_1$, $x_2$ with 
\begin{equation*}
T \leq x_1, x_2 \leq T\exp(4 \log^{20/21}{T})
\end{equation*}
such that 
\begin{align*}
\sum_{p}\varepsilon(k; p, l, 1) \log{p}\cdot  \exp\bigg( -\frac{1}{r_1}\log^2\Big(\frac p x_1 \Big)\bigg) &\geq T^{\beta}\exp \Big(-(1+\gamma^2)\log^{5/7}{T} \Big)\\
\sum_{p}\varepsilon(k; p, l, 1) \log{p}\cdot  \exp\bigg( -\frac{1}{r_2}\log^2\Big(\frac p x_2 \Big)\bigg) &\leq -T^{\beta}\exp \Big(-(1+\gamma^2)\log^{5/7}{T} \Big)
\end{align*}
Again with the contribution of primes $p$ with $p>T \exp(\log^{41/42} T)$ and $p<T \exp(-\log^{41/42}T)$ is $o(\sqrt T)$;
\end{thm}

\begin{thm}[\cite{1964.Turan_2} Theorem V] 
Under the conditions \eqref{IIb 5.4} and \eqref{IIb 5.5} there exist $U_1, U_2, U_3$ and $U_4$ with 
\begin{align*}
T\exp(-5\log^{20/21}T)\leq U_1< U_2\leq T\exp(5\log^{20/21}T)\\
T\exp(-5\log^{20/21}T)\leq U_3< U_4\leq T\exp(5\log^{20/21}T)
\end{align*}
such that
\begin{align*}
\sum_{\substack{U_1\leq p \leq U_2\\ p\equiv l \mod k}}{1}-\sum_{\substack{U_1\leq p \leq U_2\\ p\equiv 1 \mod k}}{1}&> T^\beta \exp\Big( (2+\gamma^2)\log^{5/7}T \Big) \\
\sum_{\substack{U_1\leq p \leq U_2\\ p\equiv l \mod k}}{1}-\sum_{\substack{U_1\leq p \leq U_2\\ p\equiv 1 \mod k}}{1}&<- T^\beta \exp\Big( (2+\gamma^2)\log^{5/7}T \Big)
\end{align*}
\end{thm}

Now Theorem~\ref{IIb 3} is a special case of:
\begin{thm}[\cite{1964.Turan_2} Theorem VI] 
For a $k$ satisfying the HC (Conjecture~\ref{HC}) prescribe quadratic residue $l_1$ and quadratic non-residue $l_2 \mod k$ with no $L(s, \chi)$ vanishes  for $\sigma>\h$ with $\chi(l_1)\neq \chi(l_2)$, then for suitable $c_{34}, c_{35}, c_{36}$ and
\begin{equation*}
r_0=c_{34}\frac{\log k}{E(k)^2}
\end{equation*}
the inequalities
\begin{equation*}
\sum_p\varepsilon(k; p, l_2, l_1)\log{p}\exp\Big(-\frac{1}{r}\log^2\frac{p}{x} \Big) > c_{35}\sqrt{x}
\end{equation*}
holds whenever
\begin{equation*}
r_0\leq r \leq \log x
\end{equation*}
and 
\begin{equation*}
x>c_{36}k^{50}
\end{equation*}
\end{thm}

We now present the case when
\begin{equation}\label{IIIb 2.1}
l_1=1,\;\;\;\;l_2=l={\rm quadratic \; residue\; mod\;}k
\end{equation}
\begin{thm}[\cite{1965.Turan_3} Theorem I]
For $k$ satisfies the HC (Conjecture~\ref{HC}) and in case of equation \eqref{IIIb 2.1} and for
\begin{equation}\label{IIIb 2.2}
T>\max\bigg(c_{37}, \exp\big(4\exp(3k)\big), \exp\Big(\frac{(20\pi)^6}{E(k)^6} \Big)\bigg)
\end{equation}
there exist $x_1$, $x_2$ in the interval 
\begin{equation*}
\Big( T\exp\big(-(\log T)^{5/6}\big), T\exp\big((\log{T})^{11/15}\big) \Big)
\end{equation*}
such that for suitable 
\begin{equation*}
(2\log T)^{2/3}\leq \nu_1, \nu_2 \leq (2\log T)^{2/3}+(2\log T)^{2/5}
\end{equation*}
both the inequalities
\begin{align*}
\sum_p\varepsilon(k; p, l_2, l_1)\log{p}\exp\Big(-\frac{1}{\nu_1}\log^2\frac{p}{x_1} \Big) & > \sqrt{T}\exp\big(-c_{37}\log^{5/6}T\big)\\
\sum_p\varepsilon(k; p, l_2, l_1)\log{p}\exp\Big(-\frac{1}{\nu_2}\log^2\frac{p}{x_2} \Big) &<- \sqrt{T}\exp\big(-c_{37}\log^{5/6}T\big)
\end{align*}
hold.
\end{thm}
This is a special case of:
\begin{thm}[\cite{1965.Turan_3} Theorem II]
In case \eqref{IIIb 2.1} for $k$ holding the HC (Conjecture~\ref{HC}), if $\varrho = \beta+i\gamma$ is a zero of an $L(s, \chi)$ \mod k with
\begin{equation*}
\beta\geq \h, \;\; \gamma>0 \;\;, \chi(l)\neq 1
\end{equation*}
there exist for
\begin{equation*}
T>\max \bigg(c_{38}, \exp\Big( 4\exp(3k )\Big), \exp\Big( \frac{(20\pi)^6}{E(k)^6} \Big), \exp\Big(\exp(10|\varrho|)\Big) \bigg)
\end{equation*}
$x_1, x_2$ in the interval:
\begin{equation*}
T\exp\Big(-(\log T)^{5/6}\Big) < x_1,  x_2 <  T\exp\Big((\log T)^{11/15}\Big)
\end{equation*}
such that both the inequalities
\begin{align*}
\sum_p\varepsilon(k; p, l_2, l_1)\log{p}\exp\Big(-\frac{1}{r_1}\log^2\frac{p}{x_1} \Big) & > T^{\beta}\exp\big(-c_{39}\log^{5/6}T\big)\\
\sum_p\varepsilon(k; p, l_2, l_1)\log{p}\exp\Big(-\frac{1}{r_2}\log^2\frac{p}{x_2} \Big) &<- T^{\beta}\exp\big(-c_{39}\log^{5/6}T\big)
\end{align*}
hold.
\end{thm}

\begin{thm}[\cite{1965.Turan_3} Theorem III]
For a $k$ satisfies HC (Conjecture~\ref{HC}) and in the case \eqref{IIIb 2.1} for $T$'s satisfying \eqref{IIIb 2.2} there exist numbers $U_1, U_2, U_3$ and $U_4$ with
\begin{align*}
T\exp\Big(-(\log^{6/7} T) \Big) \leq U_1 &< U_2 \leq T\exp\Big( (\log^{6/7} T) \Big)\\
T\exp\Big(-(\log^{6/7} T) \Big) \leq U_3 &< U_4 \leq T\exp\Big( (\log^{6/7} T) \Big)
\end{align*}
such that
\begin{align*}
\sum_{U_1 \leq p \leq U_2} \varepsilon(k; p, 1, l)&>\sqrt{T}\exp\bigg(-c_{40}\log^{5/6}T \bigg)\\
\sum_{U_3 \leq p \leq U_4} \varepsilon(k; p, 1, l)&<-\sqrt{T}\exp\bigg(-c_{40}\log^{5/6}T \bigg)
\end{align*}
\end{thm}

Now passing to more general cases, as we showed that more primes $\equiv l_1\mod k$ than $\equiv l_2 \mod k$ if and only if $l_1$ is an quadratic non-residue and $l_2$ is  quadratic residue $\mod k$

Let $k$ satisfy the HC (Conjecture~\ref{HC}), compare the residue classes 
\[
\equiv l_1 \mod k\;\; {\rm and}\;\; \equiv l_2 \mod k
\]
when $l_1$ and $l_2$ are both quadratic non-residues, but with more conditions: we need an $\eta$ and a small positive constant $c_{41}$ with the condition
\begin{equation}\label{IVb 1.3}
0<\eta<\min\Big(c_{41}, \Big(\frac{E(k)}{6\pi} \Big)^2 \Big)
\end{equation}
the non-vanishing of all $L(s, \chi)$ functions $\mod k$ for
\begin{equation}\label{IVb 1.4}
\sigma>\h, \;\; |t|\leq \frac{2}{\sqrt \eta}
\end{equation}
And we assume without the loss of generality that 
\begin{equation}\label{IVb 1.5}
E(k)\leq \frac{1}{k^{15}}
\end{equation}

\begin{thm}[\cite{1965.Turan_4} Theorem I]\label{IVb I}
If for $k> c_{42}$ with $c_{42}$ large and satisfying the above conditions, then for
\begin{equation*}
T>\max\bigg(c_{43}, \exp\Big( \frac{2}{\eta^4}\exp\big( \frac{1}{4}k^{10}\big)\Big) \bigg)
\end{equation*}
and for quadratic non-residue $l_1$ and $l_2$ there are $x_1, x_2, \nu_1$ and $\nu_2$ with

\begin{equation*}
T^{1-\sqrt \eta} \leq x_1, x_2 \leq T \exp(\log^{3/4}T)
\end{equation*}
and
\begin{equation*}
2\eta \log T\leq \nu_1, \nu_2 \leq 2\eta \log T +\sqrt{\log  T}
\end{equation*}
so that 
\[
\sum_{p \equiv l_1 \mod k} \log{p}\exp\Big(-\frac{1}{\nu_1}\log^2\frac{p}{x_1} \Big)-
\sum_{p \equiv l_2 \mod k} \log{p}\exp\Big(-\frac{1}{\nu_1}\log^2\frac{p}{x_1} \Big) > T^{\h -4\sqrt \eta} 
\]
\end{thm}

\begin{thm}[\cite{1965.Turan_4} Theorem II]
Under the assumptions of the previous Theorem~\ref{IVb I} there are $\mu_1, \mu_2, \mu_3, \mu_4$ with
\begin{align*}
T^{1 -4\sqrt \eta} \leq \mu_1&< \mu_2\leq T^{1 + 4\sqrt \eta} \\
T^{1 -4\sqrt \eta} \leq \mu_3&< \mu_4\leq T^{1 + 4\sqrt \eta}
\end{align*}
so that
\begin{align*}
\sum_{\substack{p\equiv l_1 \mod k\\ \mu_1 \leq p \leq \mu_2 }}1 - \sum_{\substack{p\equiv l_2 \mod k\\ \mu_1 \leq p \leq \mu_2 }} 1&> T^{\h - 5\sqrt \eta} \\
\sum_{\substack{p\equiv l_1 \mod k\\ \mu_3 \leq p \leq \mu_4 }}1 - \sum_{\substack{p\equiv l_2 \mod k\\ \mu_3 \leq p \leq \mu_4 }} 1&<- T^{\h - 5\sqrt \eta} 
\end{align*}
\end{thm}

\begin{thm}[\cite{1965.Turan_5} Theorem]
 If for a $\delta$ with $0<\delta < \frac{1}{10}$  and  for 
 \begin{equation*}
 k>\max\big( c_{44}, \exp(\delta^{-20})\big)
 \end{equation*}
 where no $L(s,\chi)$ with $\chi(l)\neq 1$, mod $k$, vanishes for
 \begin{equation*}
 |s-1|\leq \h +4\delta
 \end{equation*}
 then if
 \begin{equation*}
 a>\max\big( c_{45}, \exp(k\log^3k)\big)
 \end{equation*}
 and
 \begin{equation*}
 b=\exp\big(\log^2a \cdot(\log_2 a)^3  \big)
 \end{equation*}
 then we have $x_1, x_2$ where
 $$
 a\leq x_1, x_2<b $$
  such that
 \begin{align*}
\sum_{\substack{n \leq x_1\\ n \equiv 1\mod k}} \Lambda(n)-  \sum_{\substack{n\leq x_1\\ n\equiv l\mod k}} \Lambda(n) &\geq x_1^{\h-4\delta}\\
\sum_{\substack{n \leq x_2\\ n \equiv 1\mod k}} \Lambda(n)-  \sum_{\substack{n\leq x_2\\ n\equiv l\mod k}} \Lambda(n) &\leq -x_2^{\h-4\delta}
 \end{align*}
 \end{thm}

We return to ``modified Abelian means'', i.e. to compare between the number of primes belonging to progression $\equiv l_1 \mod k$ and $\equiv l_2 \mod k$, where both $l_1$ and $l_2$ are quadratic residues $\mod k$ 
\begin{thm}[\cite{1966.Turan_6} Theorem I]\label{VI 1}
For $l_1$, $l_2$ with $(l_1, k)=(l_2, k)=1$, $l_1\not\equiv l_2 \mod k$ are both quadratic residues $\mod k$, and conditions \eqref{V 1.3}, \eqref{IVb 1.3}, \eqref{IVb 1.4} and \eqref{IVb 1.5} hold, then for every
\begin{equation*}
T>e_2(\eta^{-3})
\end{equation*}
there are $x_1, x_2$ and $\nu_1, \nu_2$ with
\begin{align*}
T^{1-\sqrt\eta}\leq x_1,& x_2 \leq T\log T\\
2\eta \log T\leq \nu_1, & \nu_2\leq 2\eta \log T+\log_2 T
\end{align*}
such that:
\begin{align*}
\sum_{p \equiv l_1 \mod k} \log{p}\exp\Big(-\frac{1}{\nu_1}\log^2\frac{p}{x_1} \Big)-
\sum_{p \equiv l_2 \mod k} \log{p}\exp\Big(-\frac{1}{\nu_1}\log^2\frac{p}{x_1} \Big) &> T^{\h -2\sqrt \eta} \\
\sum_{p \equiv l_1 \mod k} \log{p}\exp\Big(-\frac{1}{\nu_2}\log^2\frac{p}{x_2} \Big)-
\sum_{p \equiv l_2 \mod k} \log{p}\exp\Big(-\frac{1}{\nu_2}\log^2\frac{p}{x_2} \Big) &< -T^{\h -2\sqrt \eta}
\end{align*}
hold.
\end{thm}

Analogously in short intervals we have:
\begin{thm}[\cite{1966.Turan_6} Theorem II]
Under the assumptions of the previous Theorem~\ref{VI 1} there are $\mu_1, \mu_2, \mu_3, \mu_4$ with
\begin{align*}
T^{1 -4\sqrt \eta} \leq \mu_1&< \mu_2\leq T^{1 + 4\sqrt \eta} \\
T^{1 -4\sqrt \eta} \leq \mu_3&< \mu_4\leq T^{1 + 4\sqrt \eta}
\end{align*}
so that
\begin{align*}
\sum_{\substack{p\equiv l_1 \mod k\\ \mu_1 \leq p \leq \mu_2 }}1 - \sum_{\substack{p\equiv l_2 \mod k\\ \mu_1 \leq p \leq \mu_2 }} 1&> T^{\h - 3\sqrt \eta}, \\
\sum_{\substack{p\equiv l_1 \mod k\\ \mu_3 \leq p \leq \mu_4 }}1 - \sum_{\substack{p\equiv l_2 \mod k\\ \mu_3 \leq p \leq \mu_4 }} 1&<- T^{\h - 3\sqrt \eta}. 
\end{align*}
\end{thm}

\begin{thm}[\cite{1972.Turan_7} Theorem]
There exist numbers $U_1, U_2, U_3 , U_4$ for $T>c_{46}$ with
\begin{align*}
\log_3 T \leq U_2 \exp(-\log^{15/16}{U_2})\leq U_1 < U_2 \leq T, \\
\log_3 T \leq U_4 \exp(-\log^{15/16}{U_4})\leq U_3 < U_4 \leq T, 
\end{align*}
such that
\begin{align*}
\sum_{\substack{U_1<p<U_2\\ p\equiv 1 \mod 4}}\log p - \sum_{\substack{U_1<p<U_2\\ p\equiv 3 \mod 4}}\log p& >\sqrt{U_2},\\
\sum_{\substack{U_3<p<U_4\\ p\equiv 1 \mod 4}}\log p - \sum_{\substack{U_3<p<U_4\\ p\equiv 3 \mod 4}}\log p &<-\sqrt{U_4},
\end{align*}
\end{thm}
providing insights to Problem~\ref{P:S-L}.
\section{More Chebyshev Type Assertions}
Several authors, remarkably J. Besenfelder~\cite{1979.Besenfelder_1}~\cite{1980.Besenfelder_2} and H. Bentz~\cite{1982.Bentz} proved a few unconditional theorems in the flavour of Chebyshev's assertion~\eqref{eq.cheb}
\begin{thm}[\cite{1979.Besenfelder_1} Theorem]
\[
\lim_{x \to \infty}\sum_{p}(-1)^{(p-1)/2}\log p\cdot p^{-1/2}\exp\big({-{\log^2p}/{4x}}\big)=-\infty
\]
\end{thm}
which is a special care of:
\begin{thm}[\cite{1980.Besenfelder_2} Theorem]
\[
\lim_{x \to \infty}\sum_p (-1)^{{(p-1)}/{2}}\log{p}\cdot p^{-\alpha} \exp\big(- {\log^2 p}/{4x}\big)=-\infty\;\;\;\; {\rm for}\;\;\;\; 0 \leq \alpha \leq \h
\]
Where the magnitude of divergence for $\alpha = \h$ is given by $\h\sqrt{\pi y}$ and for $0\leq \alpha<\h $ is given by $\sqrt{\pi y}e^{\frac{y}{4}(1-2\alpha)}$
\end{thm}
Further, H. Bentz proves~\cite{1982.Bentz}
\begin{thm}[\cite{1982.Bentz} Theorem 1]
Unconditionally,
\[
\lim_{x \to \infty}\sum_{p}(-1)^{(p-2)/2}\log{p}\cdot{p^{-1/2}}\exp\big(-\log^2p/x \big)=-\infty 
\]
The magnitude of divergence is given by $\frac{1}{4}\sqrt{\pi x}+O(1)$.
\end{thm}
which generalizes to: 
\begin{thm}[\cite{1982.Bentz} Theorem 2]
\[
\lim_{x \to \infty}\sum_p (-1)^{(p-1)/2}\log p\cdot{p^{-\alpha}}\exp(-\log^2p/x)=-\infty\;\;\;\;{\rm for\; all}\;\;\;\; 0\leq \alpha < \h
\]
The magnitude of divergence is given by 
$
\sim\h \sqrt{\pi x} \exp\Big(\frac{x}{16}(1-2\alpha)^2\Big)
$
\end{thm}
\begin{dfn}\label{chi_3}
If we define
\[
\chi_3(m)=\begin{cases} 1 \;\;\;\; &{\rm if} \;\;\; m\equiv 1 \mod 3,\\
                       -1 \;\;\;\; &{\rm if} \;\;\; m\equiv 2 \mod 3,\\
                        0 \;\;\;\; &{\rm if} \;\;\; m\equiv 0 \mod 3
                        \end{cases}
                        \]
                        i.e. taking $k=3$ and thus $\varphi(k)=2$
\end{dfn}
We have:
\begin{thm}[\cite{1982.Bentz} Theorem 3]
Let $\chi_3$ be given as in Definition~\ref{chi_3}, then
\[
\lim_{x\to \infty}\sum_{p}\chi_3(p)\log p\cdot p^{-1/2}\exp(-\log^2p/x)=-\infty
\]
The order of magnitude of divergence is given by
$\frac{1}{4}\sqrt{\pi x}+O(1)
$\end{thm}
\begin{thm}[\cite{1982.Bentz} Theorem 4]
Let $\chi_3$ be as in Definition~\ref{chi_3}, then
\[
\lim_{x\to \infty}\sum_p\chi_3(p){\log p}\cdot{p^\alpha}\exp(-\log^2p/x)=-\infty\;\;\;\;{\rm{for}\;} 0 \leq \alpha < \h.
\]
The order of magnitude of divergence is given by
$\sim \h\sqrt{\pi x}\exp\big({\frac{x}{16}(1-2\alpha)^2}\big)
$\end{thm}
Now for a higher moduli, in the cases $k=8$ and $5$ so $\varphi(k)=4$, H. Bentz and J. Pintz prove:
\begin{thm}[\cite{1982.Bentz} Theorem 5]
\begin{align*}
\lim_{x \to \infty}\sum_{p \equiv 1 \mod 8}\log p\cdot p^{-1/2}\exp(-\log^2p/x)-\sum_{p \equiv 3 \mod 8}\log p\cdot p^{-1/2}\exp(-\log^2p/x)=- \infty\\
\lim_{x \to \infty}\sum_{p \equiv 1 \mod 8}\log p\cdot p^{-1/2}\exp(-\log^2p/x)-\sum_{p \equiv 5 \mod 8}\log p\cdot p^{-1/2}\exp(-\log^2p/x)=- \infty\\
\lim_{x \to \infty}\sum_{p \equiv 1 \mod 8}\log p\cdot p^{-1/2}\exp(-\log^2p/x)-\sum_{p \equiv 7 \mod 8}\log p\cdot p^{-1/2}\exp(-\log^2p/x)=- \infty\\
\end{align*}
with the order of magnitude of divergence being $-\frac{1}{4}\sqrt{\pi y}+O(1)$, respectively.
\end{thm}
\begin{thm}[\cite{1982.Bentz} Theorem 6]
\begin{align*}
\lim_{x \to \infty}\bigg\{\sum_{p\equiv 3 \mod 8}-\sum_{p\equiv 5 \mod 8}\bigg\}\log p\cdot p^{-1/2}  \exp(-\log^2p/x)=O(1)\\
\lim_{x \to \infty}\bigg\{\sum_{p\equiv 3 \mod 8}-\sum_{p\equiv 7 \mod 8}\bigg\}\log p\cdot p^{-1/2}  \exp(-\log^2p/ x)=O(1)\\
\lim_{x \to \infty}\bigg\{\sum_{p\equiv 5 \mod 8}-\sum_{p\equiv 7 \mod 8}\bigg\}\log p\cdot p^{-1/2} \exp(-\log^2p/x)=O(1)
\end{align*}
\end{thm}

\begin{thm}[\cite{1982.Bentz} Theorem 7]
For at least one of the two classes $2 \mod 5$, $3 \mod 5$, we have
\[
\lim_{x\to \infty}\Bigg(\bigg(\sum_{\substack{p \equiv 2 \mod 5\\{\rm or}\; p \equiv 3 \mod 5}}-\sum_{p \equiv 4 \mod 5}\bigg)\log p \cdot p^{-1/2}\exp(-\log^2 p/4x )\Bigg) =+\infty
\]
\end{thm}

and when dealing with quadratic residues and distribution of primes, H. Bentz~\cite{1980.Bentz} assumes the following two conjectures:
\begin{cnj}[$\bf{R}_2$]\label{R_2}
The domain $\sigma > \h$, $|t|\leq 1$ is zero free and there is NO zero at $s=\h$ for Dirichlet $L$-function.
\end{cnj}

\begin{cnj}[$\bf{H}_2$]\label{H_2}
All zeros $\varrho:=\beta+i \gamma$ satisfy the inequality 
\[
\beta^2-\gamma^2<\frac{1}{4}
\]
\end{cnj}
and he shows
\begin{thm}[\cite{1980.Bentz} Theorem 1]
If $l_1$ is a quadratic residue, $l_2$ a non-residue mod $k$ and ${\bf R_2}$ (Conjecture~\ref{R_2}) or even ${\bf H_2}$ (Conjecture~\ref{H_2}) valid for $L$-function mod $k$, then 
\[
\lim_{x \to \infty}\sum_p\varepsilon (k; p, l_1, l_2)\log p\cdot\exp(-{\log^2 p}/{x} ) =-\infty
\]
\end{thm}

\begin{thm}[\cite{1980.Bentz} Theorem 2]
If $l_1$ a quadratic residue, $l_2$ a non-residue mod $k$, then
\[
\lim_{x \to \infty} \sum_{p}\varepsilon(k; p, l_1, l_2)\log p\cdot\exp(-{\log^2 p}/{x}) = - \infty
\] 
holds for all $k< 25$.
\end{thm}

\begin{thm}[\cite{1980.Bentz} Theorem 3]\label{1980.Bentz Thm 3}
If ${\bf{ R_2}}$ (Conjecture~\ref{R_2}) or only ${\bf{ H_2}}$ (Conjecture~\ref{H_2}) is true for all $L$-functions $\mod k$, $l_1$ is a quadratic residue, $l_2$ a non-residue mod $k$, then for $0 \leq \alpha < \h$
\[
\lim_{x \to \infty}\sum_{p}\varepsilon(k, p, l_1, l_2)\log p \cdot p^{-\alpha}\exp(-\log^2 p/x) = -\infty
\]
\end{thm}

\begin{thm}[\cite{1980.Bentz} Theorem 4]\label{1980.Bentz Thm 4}
If $l_1$ a quadratic residue, $l_2$ a non-residue $\mod k$, $k< 25$, then for $0 \leq \alpha < \h$
\[
\lim_{x \to \infty}\sum_p \varepsilon(k; p, l_1, l_2)\log p \cdot p^{-\alpha}\exp(-\log^2 p/{x}) = -\infty
\]
\end{thm}

\begin{thm}[\cite{1980.Bentz} Theorem 5]
Under the condition of Theorem~\ref{1980.Bentz Thm 3} we have
\[
\sum_n \varepsilon(k; n, l_1, l_2)\log p\cdot{p^{-\alpha}} \exp( -\log^2 p/x)\sim \frac{N(k)}{\varphi(k)}\sqrt{\pi x}\exp\Big((x/4)(1/2 -x)^2 \Big)
\]
where $N(k)$ denotes the number of solutions of $x^2 \equiv 1 \mod k$
\end{thm}

Of course the above theorem implies:
\begin{thm}[\cite{1980.Bentz} Theorem 6]
Under the conditions of Theorem~\ref{1980.Bentz Thm 4}
\[
\sum_n \varepsilon (k, n, l_1, l_2)\log p\cdot p^{-\alpha}\exp(-\log^2 p/x) \sim \frac{N(q)}{\varphi(q)}\sqrt{\pi x}\exp\Big((x/4)(1/2-x)^2 \Big)
\]
\end{thm}

\section{A Few Other Results}
Knapowski and Tur\'an also made contributions to Problem~\ref{P:L-B} for $\Delta_\pi(r; k, l_1, l_2)$ in the cases of $k=8$ and $4$:
\begin{thm}[\cite{1965.Knapowski} Theorem I]
For any $l_1 \neq l_2$ among $3, 5, 7$ and $0<\delta<c_{47}$, we have the inequality
\begin{equation*}
\max_{\delta\leq x \leq \delta ^{{1/3}}} |\Delta_\pi(r; 8, l_1, l_2)| \geq \delta^{-1/2} \exp\bigg(\frac{23\log(1/\delta)\log_3(1/\delta)}{\log_2({1/\delta})}\bigg)
\end{equation*}
\end{thm}

\begin{thm}[\cite{1965.Knapowski} Theorem II]
For $l \neq 1, k=4$ \emph{or} $8$ and $0<\delta< c_{48}$,
\begin{equation*}
\max_{\delta\leq x \leq \delta ^{{1/3}}} |\Delta_\pi(r; k, 1, l)| \geq \delta^{-1/2} \exp\bigg(\frac{23\log(1/\delta)\log_3(1/\delta)}{\log_2({1/\delta})}\bigg)
\end{equation*}
\end{thm}
\vskip12pt
In his paper~\cite{1971.Stark}, H. Starks studies the asymptotic behaviours of $\varphi(k)\pi(x, k, a)- \varphi(K)\pi(x, K, A)$: If $\chi$ and $X$ are characters mod $k$ and $K$ respectively, and $\chi_0$ and $X_0$ denote the principle characters, whereas $\chi_R$ and $X_R$ denote the real characters, he defines:
\begin{dfn}
\begin{equation}\label{r}
r:=r(k, a; K, A) = \sum_{X_R}{X_R(A)}- \sum_{\chi_R}{\chi_R(A)}
\end{equation}
\begin{align*}
A_T(u)&:= A_T(u, k, a, K, A)\\
           &= \sum_{X\neq X_0}\sum_{\substack{\varrho_X \\ \beta_X > 0, |\gamma_X |< T}} \frac{\overline X(A)}{\varrho_X}\exp{(\varrho_X-\h)u} - 
                \sum_{\chi \neq \chi_0} \sum_{\substack{\varrho_\chi \\ \beta_\chi > 0, |\gamma_\chi |< T}} \frac{\bar\chi(A)}{\varrho_\chi}\exp{(\varrho_\chi-\h)u} \\
A_T^*(u)&:= A_T^*(u, k, a, K, A)\\
           &= r+ \sum_{X\neq X_0}\sum_{\substack{\varrho_X \\ \beta_X > 0, |\gamma_X |< T}} \frac{X(A)}{\varrho_X}\exp{(\varrho_X-\h)u} - 
                \sum_{\chi \neq \chi_0} \sum_{\substack{\varrho_\chi \\ \beta_\chi > 0, |\gamma_\chi |< T}} \frac{\chi(A)}{\varrho_\chi}\exp{(\varrho_\chi-\h)u}
\end{align*}
so the relation between them is simply
\begin{equation*}
A_T^*(u)= r+\frac{1}{T}\int_0^T A_t(u) \, dt
\end{equation*}
further, whenever the limit exists, define
$A_\infty(u):=A_\infty(u; k; a; K, A) = \lim_{T\to \infty}A_T(u; k, a; K, A)$
\end{dfn}
\begin{thm}[\cite{1971.Stark} Theorem 1]
Under GRH (Conjecture~\ref{GRH}), for any $T>0$ and any $u$,
\[
\limsup_{y \to \infty}\frac{\varphi(k)\pi(x, k, a)- \varphi(K)\pi(x, K, A)}{\sqrt y/\log y}\geq A_T^*(u)
\]
\end{thm}

\begin{thm}[\cite{1971.Stark} Theorem 2]
Again assuming GRH (Conjecture~\ref{GRH}) 
\begin{enumerate}
\item\label{1971.1} If $r(k, a; K, A) = 0$, then there is a constant $c>0$ such that
\begin{equation*}
\limsup_{y \to \infty}\frac{\varphi(k)\pi(x, k, a)- \varphi(K)\pi(x, K, A)}{{{\sqrt{y}}/{\log y}}}\geq c
\end{equation*}
\item If $r(k, a; K, A) > 0$, then the result of (\ref{1971.1}) is true with $c=r$, $r$ in equation~\eqref{r}.
\end{enumerate}
\end{thm}
\vskip12pt
On the sign changes of $\pi(x; q, 1)-\pi(x; q, a)$, J.-C. Schlage-Puchta engenders:
\begin{thm}[\cite{2004.Schlage-Puchta} Theorem 1]
When $q$ s a natural number, we define $q^+:=\max\big(q, \exp(1260)\big)$, and assuming GRH (Conjecture~\ref{GRH}). Let $M(q)$ be the number of solution of the congruence $x^2 \equiv 1 \mod q$. Then there exists an $x$ with$x< e_2\big((q^+)^{170}+e^{18M(q)} \big)$ such that $\pi(x; q, 1)> \pi(x; q, a)$ for all $a \not\equiv 1 \mod q$.  Moreover, let $V(x)$ demote the number of sign changes of $\pi(t; q, 1)- \max_{a \not\equiv 1 \mod q}\pi(t; q, a)$ in the range $2 \leq t \leq q$, then
\[
V(x)>\frac{\log x}{\exp \big( (q^+)^{170}+e^{18M(q)}\big)}-1
\]
\end{thm}


\section{Modern Developments on the Racing Problems}
Several authors had made progresses on the Shank-R\'enyi Racing Problems (Problem~\ref{P:P-R} described in Section~\ref{intro section} and their variations), notably early on by Kaczorowski~\cite{1993.Kaczorowski} as he proposed: 
\begin{cnj}[Strong Race Hypothesis]~\label{SRH}
For each permutation $a_1, a_2, \hdots, a_{\varphi(k)}$ of the reduced set of residue classes mod $k$ the set of integers $m$ with 
\[
\pi(m, k, a_1)<\pi(m, k, a_2)<\dots<\pi(m, k, a_{\varphi(k)})
\]
has positive ``lower density''.
\end{cnj}
\begin{thm}[\cite{1993.Kaczorowski} Theorem 1]
Under GRH (Conjecture~\ref{GRH}) for Dirichlet $L$-functions mod $k$, $k\geq 3$. There exists infinitely many integers $m$ with 
\[
\pi(m, k, 1)>\max_{a \not\equiv 1 \mod k}\pi(m, k, a)
\] 
Moreover, the set of $m$'s satisfying the inequality has positive density.\\
Same statement holds true for $m$ satisfying 
\[
\pi(m, k, 1)<\min_{a \not\equiv 1 \mod k}\pi(m, k, a)
\]
\end{thm}
which is an immediate consequence of:
\begin{thm}[\cite{1993.Kaczorowski} Theorem 2]\label{1993.Kaczorowski 2}
Under GRH (Conjecture~\ref{GRH}) for $L$-functions mod $k$, $k\geq 3$, and let $u$ denote an arbitrary non-negative real number.  Then there exist constants $c_{49}=c_{49}(u)>0$ and $c_{50}=c_{50}(u)>1$ only depending on $u$, such that for every $T\geq 1$
\begin{align*}
\#\Big\{T\leq m \leq c_{50}T: \psi(m, k, 1)\geq \max_{a\not\equiv 1 \mod k}\psi(m, k ,a)+u\sqrt{m}\Big\}\geq& c_{49}T\\
\#\Big\{T\leq m \leq c_{50}T: \pi(m, k, 1)\geq \max_{a\not\equiv 1 \mod k}\pi(m, k, a)+u\frac{\sqrt{m}}{\log m}\Big\}\geq& c_{49}T\\
\#\Big\{T\leq m \leq c_{50}T: \psi(m, k, 1)\leq \max_{a\not\equiv 1 \mod k}\psi(m, k ,a)-u\sqrt{m}\Big\}\geq& c_{49}T\\
\#\Big\{T\leq m \leq c_{50}T: \pi(m, k, 1)\leq \max_{a\not\equiv 1 \mod k}\pi(m, k ,a)-u\frac{\sqrt{m}}{\log m}\Big\}\geq& c_{49}T
\end{align*} 
\end{thm}

Kaczorowski also made some progress on the racing problem~\ref{P:P-R}, with $k=5$ for $\psi(m, 5, a_i)$
\begin{thm}[\cite{1995.Kaczorowski} Theorem 1]
Assuming GRH (Conjecture~\ref{GRH}) for modoluo 5. Then for every permutation $(a_1, a_2, a_3, a_4)$ of the sequence $(1, 2, 3, 4)$ the set of $m$'s satisfying
\[
\psi(m, 5, a_1) > \psi(m, 5, a_2) > \psi(m, 5, a_3)> \psi(m, 5, a_4)
\]
has positive density.
\end{thm}

\begin{thm}[\cite{1995.Kaczorowski} Theorem 2]
Assuming GRH (Conjecture~\ref{GRH}) with $L(s, \chi)$ mod $5$, there exist three positive constants $c_{51}, c_{52}, c_{53}$ such that for every permutation $(a_1, a_2, a_3, a_4)$ of the sequence $(1, 2, 3, 4)$ and for arbitrary $T\geq 1$ we have
\[
\#\{T\leq m \leq c_{51}T: \psi(m, 5, a_1) > \hdots> \psi(m, 5, a_4),
\min_{\substack{i \neq j \\ 1\leq i,j \leq 4}}|\delta_\psi(m; 5, i, j)|\geq c_{52}\sqrt{m}\}\geq c_{53} T
\]
\end{thm}
\vskip12pt
He employs the $k$-functions bearing his name in~\cite{1996.Kaczorowski}, with
\begin{dfn} For $q>1$ a natural number, let
\[
m(q):=\begin{cases}
\h \;\;\; {\rm if} \; 2 \parallel q\\
2 \;\;\; {\rm if} \;\;8 \,|\,q \\
1 \;\;\; {\rm otherwise}
\end{cases}
\]
\[
N_q:=\frac{1}{\varphi(q)}m(q)2^{\omega(q)}
\]
where $\omega(q)$ denotes the number of distinct prime divisors of $q$.\\
Also define
\[
p^{\nu_p(q)} \| q, \;\;\;\; q_p:=qp^{-\nu_p(q)}, \;\;\;\; g_{p, q}:=\ord \;p \mod {q_p}
\]
Now let $(a, q)=1$, and denote by $\bar{a}$ the inverse of $a \mod q$: $a\bar{a}\equiv 1 \mod q$. \\
 Moreover, he put
\begin{align*}
\varrho(q, a)&:=\begin{cases} 1 \;\;\;{\rm if}\; $a$\; {\rm is\; a\; quadratic\; residue} $\mod q$\\ 
0 \;\;\; {\rm otherwise}
\end{cases}\\
\lambda(q, a)&:= \sum_{\substack{p^\alpha \|q \\ a\equiv 1 \mod {q_p}}}\frac{\log p}{p^{\alpha-1}(p-1)}+\sum_{\substack{p^\alpha|q, \alpha< \nu_p(q) \\ a \equiv 1 \mod {qp^{-\alpha}}}}\frac{\log p}{p^\alpha}\\
\delta(q, a)&:= \begin{cases} 1 \;\;\; {\rm if}\; a\equiv -1 \mod k\\
0 \;\;\;{\rm otherwise} \end{cases}
\end{align*}
\end{dfn}

Suppose $p$ a prime and that $a \mod k_p$ belongs to the cyclic multiplicity group generated by $p \mod {q_p}$.  Then denote by $l_p(a)$ the natural number uniquely determined by:
\[
1 \leq l_p(a)\leq g_{q, p}, \;\;\;\;\;\;\;\;\;p^{l_p(a)}\equiv a \mod {q_p}
\]
then set
\[
\alpha(q, a):=\sum_{p|q}\frac{\log p}{\varphi(p^{\nu_p(q)})p^{l_p(a)}}\Bigg(1-\frac{1}{p^{g_{q, p}}} \Bigg)^{-1}
\]
if there are no such primes $p$ we put $\alpha(q, a)=0$.
\begin{rmk}
An easy consequence of Dirichlet's prime number theorem is that for every $a$ to $q$ there exists a constant $b(q, a)$ such that
\[
\sum_{\substack{n\leq x \\ n\equiv a \mod q}}\frac{\Lambda(n)}{n}=\frac{1}{\varphi(q)}\log{x}+b(q, a)+o(1)
\]
as $x$ tends to infinity, where $b(q, a)$ is called the Dirichlet-Euler constant.
\end{rmk}
Finally, Kaczorowski defines the following quantities:
\begin{align*}
r^+(q, a):&=\alpha(q, a)+b(q, a)+\h \delta(q, a)\log 2 +\lambda(q, a)\\
r^-(q, a):&=\alpha(q, a)+b(q, a)+\h \delta(q, a)\log 2\\
R^+(q, a):&=r^+(q, a)-\varrho(q, a)N_q\\
R^-(q, a):&=r^-(q, a)-\varrho(q, a)N_q
\end{align*}
so he is able to prove:
\begin{thm}[\cite{1996.Kaczorowski} Theorem]
Let $k\geq 5$, $k\neq 6$ be an integer and assume the GRH (Conjecture~\ref{GRH}) $\mod k$.\\
Define permutations:
\begin{align*}
(a_2, a_3, \hdots, a_{\varphi(k)}), \;\;\; (b_2, b_3, \hdots, b_{\varphi(k)})\\
(c_2, c_3, \hdots, c_{\varphi(k)}), \;\;\; (d_2, d_3, \hdots, d_{\varphi(k)})
\end{align*}
of the set of residue classes
\[
a \mod k, \;\;\;\;\; (a, k)=1\;\;\; a\not\equiv 1 \mod k
\]
so that the following inequalities hold:
\begin{align*}
R^+(k, \bar{a}_2)&>R^+(k, \bar{a}_3)>\hdots>R^+(k, \bar{a}_{\varphi(k)}) \\
R^-(k, \bar{b}_2)&>R^-(k, \bar{b}_3)>\hdots>R^-(k, \bar{b}_{\varphi(k)})\\
r^+(k, \bar{c}_2)&>r^+(k, \bar{c}_3)>\hdots>r^+(k, \bar{c}_{\varphi(k)})\\
r^-(k, \bar{d}_2)&>r^-(k, \bar{d}_3)>\hdots>r^-(k, \bar{d}_{\varphi(k)})
\end{align*}
Then there exists a positive constant $b_0$ such that each of the sets of natural numbers each set of natural numbers
\begin{align*}
\{m\in \N: &\pi(m; k, a_2)>\hdots > \pi(m; k, a_{\varphi(k)})>\pi(m; k, 1), \\ 
&\min_{a\not\equiv b \mod k, (ab, k)=1}|\pi(m; k, a)-\pi(m; k, b)|>b_0\sqrt{m}/\log m\}\\
\{m\in \N: &\pi(m; k, 1)> \pi(m; k, a_2)>\hdots>\pi(m; k, a_{\varphi(k)}), \\ 
&\min_{a\not\equiv b \mod k, (ab, k)=1}|\pi(m; k, a)-\pi(m; k, b)|>b_0\sqrt{m}/\log m\}\\
\{m\in \N: &\psi(m; k, c_2)>\hdots > \pi(m; k, c_{\varphi(q)})>\pi(m; k, 1), \\ 
&\min_{a\not\equiv b \mod k, (ab, k)=1}|\psi(m; k, a)-\psi(m; k, b)|>b_0\sqrt{m}/\log m\}\\
\{m\in \N: &\psi(m; k, 1)> \psi(m; k, d_2)>\hdots>\psi(m; k, d_{\varphi(q)}), \\ 
&\min_{a\not\equiv b \mod k, (ab, k)=1}|\psi(m; k, a)-\psi(m; k, b)|>b_0\sqrt{m}/\log m\}
\end{align*}
has a positive density.
\end{thm}
\vskip12pt
In their ground-breaking paper~\cite{1994.Rubinstein}, M. Rubinstein and P. Sarnak resurrected the racing-problem (Problem~\ref{P:P-R}) and fully solved a few open problems with the assumption of some unproven conditions mentioned in Section~\ref{intro section}, namely GRH (Conjecture~\ref{GRH}) and:
\begin{cnj}
[Linear Independence hypothesis (LI)]\label{LI}
The imaginary part of the zeros of all Dirichlet $L$-functions attached to primitive characters modulo $q$ are linearly independent over $\Q$.
\end{cnj}
They employed the logarithmic density:
\begin{dfn}
\begin{align*}
\bar{\delta}(P):&=\limsup_{X \to \infty}\frac{1}{X}\int_{t\in P \cap [2, X]}\,\frac{dt}{t}\\
\underline{\delta}(P):&=\liminf_{X \to \infty}\frac{1}{X}\int_{t\in P \cap [2, X]}\,\frac{dt}{t}
\end{align*}
and set $\delta(P)=\bar\delta(P)=\underline\delta(P)$ if the above two limits are equal.
\end{dfn}
By introducing the vector-valued functions,
\begin{dfn}
\[
E_{k; a_1, a_2, \hdots, a_r}(x):=\frac{\log x}{\sqrt x}\times (\varphi(k)\pi(x; k, a_1)-\pi(x), \hdots, \varphi(k)\pi(x; k, a_r)-\pi(x))
\]
for $x \geq 2$.
\end{dfn}
they studied the existence of and tried to estimate the logarithmic density of of the set $P_{k; a_1, \hdots, a_r}$, where 
\begin{dfn}
 $P_{k; a_1, \hdots, a_r}$ is the set of real numbers $x \geq 2$ such that
\[
\pi(x; k, a_1) > \pi(x; q, a_2)> \dots > \pi(x; k, a_r)
\]
with $k\geq 3$ and $2 \leq r \leq \varphi(k)$, and denote $\A_r(k)$ the set of ordered $r$-tuples of distinct residue classes $(a_1, a_2, \hdots, a_r)$ modulo $k$ which are coprime to $k$. 
\end{dfn}
so they could the following theorems:
\begin{thm}[\cite{1994.Rubinstein} Theorem 1.1]
Under GRH (Conjecture~\ref{GRH}), $E_{k; a_1, a_2, \hdots, a_r}$ has a limiting distribution $\mu_{k; a_1, \hdots, a_r}$ on $\R^r$, i.e. 
\[
\lim_{X \to \infty}\frac{1}{X}\int_2^X f(E_{k; a_1, a_2, \hdots, a_r}(x)) \frac{dx}{x}=\int_{\R^r}f(x)\,d\mu_{k; a_1, \hdots, a_r}(x)
\]
for all bounded continuous functions $f$ on $\R^r$.
\end{thm}
\begin{rmk}
If it turns out that if the measure $\mu_{k; a_1, \hdots, a_r}$ is absolutely continuous then 
\[
\delta(P_{k; a_1, \hdots, a_r})=\mu_{k; a_1, \hdots, a_r}\big(\{x \in \R^r: x_1>\cdots>x_r \}\big)
\]
the shortcoming here is that they write this assuming only GRH (Conjecture~\ref{GRH}), they do not know that $\delta(P_{k; a_1, \hdots, a_r})$ exists.
\end{rmk}

\begin{dfn}
Since the measures $\mu$ are very localized but not compactly supported:
\begin{align*}
B'_R:&= \{x\in \R^r : |x|\geq R\}\\
B^+_R:&=\{x \in B'_R : \varepsilon(a_j)x_j>0 \}\\
B^-_R:&=-B^+_R
\end{align*}
\[{\rm {where}}\;\; \varepsilon(a)= \begin{cases} 1 &\;\;\; \rm{if}\;\; a \equiv 1 \mod k\\
                                                   -1 &\;\;\; \rm{otherwise}
                                                    \end{cases}\]
\end{dfn}

\begin{thm}[\cite{1994.Rubinstein} Theorem 1.2]
With GRH (Conjecture~\ref{GRH}), there are positive constants $c_{54}, c_{55}, c_{56}, c_{57}$ depending only on $k$ such that 
\begin{align*}
\mu_{k; a_1, \hdots, a_r}(B'_R)\leq& c_{54} \exp(-c_{55}\sqrt{R})\\
\mu_{k; a_1, \hdots, a_r}(B^{\pm}_R)\geq& c_{56} \exp(-\exp c_{57}{R})
\end{align*}
\end{thm}

H. L. Montgomery~\cite{1980.Montgomery} under RH (Conjecture~\ref{RH}) and LI (Conjecture~\ref{LI}) for $\zeta(s)$, investigated the tails of the measure $\mu_{1; 1}$, where he showed
\[
\exp\big(-c_{58}\sqrt{R}\exp(\sqrt{2\pi R}\,)\big)\leq \mu_{1;1}(B^\pm_R)\leq  \exp\big(-c_{59}\sqrt{R}\exp(\sqrt{2\pi R}\,)\big)
\]
Rubinstein and Sarnak \cite{1994.Rubinstein} under GRH (Conjecture~\ref{GRH}) and LI (Conjecture~\ref{LI}) have found an explicit formula for the Fourier transform of $\mu_{k; a_1, \hdots, a_r}$: the formula says that, for $r<\varphi(k)$, $\mu_{k; a_1, \hdots, a_r}=f(x)\,dx$ with a rapidly decreasing entire function $f$.  As a consequence, under GRH (Conjecture~\ref{GRH}) and LI (Conjecture~\ref{LI}) each $\delta(P_{k; a_1, \hdots, a_r})$ does indeed exist and is non-zero (including the case $r=\varphi(k)$).  Therefore, the solution to the racing problem~\ref{P:P-R} is conditionally affirmative.

\begin{dfn}
Define $(k; a_1, \hdots, a_r)$ to be {\bf{unbiased}}, if the density function of $\mu_{k; a_1, \hdots, a_r}$ is invariant under permutations of $(x_1, \hdots, x_r)$.
Where
\[
\delta(P_{k; a_1, \hdots a_r})=\frac{1}{r!}
\]
further define
\[
c(q, a):=-1 +\sum_{\substack{b^2 \equiv a \mod q \\ 0\leq b \leq q-1}}1
\]
\end{dfn}
\begin{thm}[\cite{1994.Rubinstein} Theorem 1.4]
Assuming GRH (Conjecture~\ref{GRH}) and LI (Conjecture~\ref{LI}) for $\chi \mod k$, $(k; a_1, \hdots a_r)$ is unbiased if and only if either $r=2$ and $c(k, a_1)=c(k, a_2)$ or $r=3$ and there exists $\rho \neq 1$ such that $\rho^3 \equiv 1 \mod k$,\; $a_2 \equiv a_1\rho \mod k$, and \; $a_3 \equiv a_1\rho^2 \mod k$.
\end{thm}

\begin{thm}[\cite{1994.Rubinstein} Theorem 1.5] Assuming GRH (Conjecture~\ref{GRH}) and LI (Conjecture~\ref{LI}) modulo $k$, for $r$ fixed,
\[
\max_{a_1, \hdots, a_r \in \mathcal{A}_q}\bigg|\delta(P_{k; a_1, \hdots, a_r})-\frac{1}{r!} \bigg| \to 0 \;\;{\rm{as}}\;\;q\to \infty 
\]
\end{thm}
\begin{thm}[\cite{1994.Rubinstein} Theorem 1.6]
Assume GRH (Conjecture~\ref{GRH}) and LI (Conjecture~\ref{LI}). Let $\tilde{\mu}_{k; N, R}$ be the limiting distribution of 
\[
\frac{E_{k; N, R}(x)}{\sqrt{\log q}}
\]
then $\tilde{\mu}_{k; N, R}$ converges in measure to the Gaussian $(2\pi)^{-1/2}\exp(-x^2/2)\,dx$ as $q \to \infty$.
\end{thm}
\vskip12pt
A. Feuerverger and G. Martin's paper~\cite{2000.Feuerverger} first presents some \emph{biased} examples using Rubinstein and Sanark's notation $\delta_{k; a_1, \hdots, a_r}$ with numerical values: for $k=8$ and $12$:
\begin{thm}[\cite{2000.Feuerverger} Theorem 1]
Assume GRH (Conjecture~\ref{GRH}) and LI (Conjecture~\ref{LI}). Then
\begin{align*}
\delta_{8; 3, 5, 7}=\delta_{8; 7, 5, 3}=0.1928013 \pm 0.000001\\
\delta_{8; 3, 7, 5}=\delta_{8; 5, 7, 3}=0.1664263 \pm 0.000001\\
\delta_{8; 5, 3, 7}=\delta_{8; 7, 3, 5}=0.1407724 \pm 0.000001
\end{align*}
and
\begin{align*}
\delta_{12; 5, 7, 11}=\delta_{12; 11, 7, 5}=0.1984521 \pm 0.000001\\
\delta_{12; 5, 11, 7}=\delta_{12; 7, 11, 5}=0.1215630 \pm 0.000001\\
\delta_{12; 7, 5, 11}=\delta_{12; 11, 5, 7}=0.1799849 \pm 0.000001
\end{align*}
where the indicated error bounds are rigorous.
\end{thm}
\begin{thm}[\cite{2000.Feuerverger} Theorem 2]
Assume GRH (Conjecture~\ref{GRH}) and LI (Conjecture~\ref{LI}), and let $k, r \geq 2$ be integers and let $a_1, \hdots, a_r$ be distinct reduced residue classes modulo $k$.
\begin{enumerate}
\item Letting $a_j^{-1}$ denote the multiplicative inverse of $a_j$ modulo $k$, we have $\delta_{k; a_1, \hdots, a_r}=\delta_{k; a_1^{-1}, \hdots, a_r^{-1}}$.
\item If $b$ is a reduced residue class modulo $k$ such that $c(k, a_j)=c(k, ba_j)$ for each $1\leq j \leq r$, then $\delta_{k; a_1, \hdots, a_r}=\delta_{k; ba_1, \hdots, ba_r}$. In particular, this holds if $b$ is a square modulo $k$.
\item If the $a_j$ are all squares modulo $k$ and $b$ is any reduced residue class modulo $k$, then $\delta_{k; a_1, \hdots, a_r}=\delta_{k; ba_1, \hdots, ba_r}$.
\item If the $a_j$  are either all squares modulo $k$ or all non-squares modulo $k$, then $\delta_{k; a_1, \hdots, a_r}=\delta_{k; a_r, \hdots, a_1}$.
\item If $b$ is a reduced residue class modulo $k$ such that $c(k, a_j) \neq c(k, ba_j)$ for each $1 \leq j \leq r$, then $\delta_{k; a_1, \hdots, a_r} = \delta_{k; ba_r, \hdots, ba_1}$. In particular, this holds if $k$ is an odd prime power or twice an odd prime power and $b$ is any non square modulo $k$.
\end{enumerate}
\end{thm}
\begin{thm}[\cite{2000.Feuerverger} Theorem 3]
Under GRH (Conjecture~\ref{GRH}) and LI (Conjecture~\ref{LI}) for $k \geq 2$ be an integer, let $N$ and $N'$ be distinct (invertible) non-squares modulo $k$, and let $S$ and $S'$ be distinct (invertible) squares $\mod k$. Then
\begin{enumerate}
\item $\delta_{k; N, N', S} > \delta_{k; S, N', N}$;
\item $\delta_{k; N, S, S'} > \delta_{k; S', S, N}$;
\item $\delta_{k; N, S, N'} > \delta_{k; N', S, N}$ if and only if $\delta_{k; N, S} > \delta_{k; N', S}$
\item $\delta_{k; S, N, S'} > \delta_{k; S', N, S}$ if and only if $\delta_{k; S, N} > \delta_{k; S', N}$
\end{enumerate}
\end{thm}

\begin{thm}[\cite{2000.Feuerverger} Theorem 4]
Assume GRH (Conjecture~\ref{GRH}) and LI (Conjecture~\ref{LI}) for $\chi \mod k$ with $k\geq 2$. Let $r \geq 2$ be an integer, and let $a_1, \hdots, a_r$ be distinct residue classes mod $k$. Then
\[
\delta_{k; a_1, \hdots, a_r}=2^{-(r-1)}\Bigg(1+\sum_{\substack{B\subset\{1, \hdots, r-1\}\\B \neq \emptyset}}\bigg(\frac{i}{\pi} \bigg)^{|B|}\times {\rm P.V.}\int\cdots\int{\hat{\varrho}_{k; a_1, \hdots, a_r(B)}\prod_{j\in B}\frac{d\eta_j}{\eta_j}} \Bigg)
\]
where $\hat{\varrho}_{k; a_1, \hdots, a_r}(B)$ borrows the notation
\[
f(B)=f(B)(\{x_j: j\in B\})=f(\theta_1, \hdots, \theta_n)
\] 
\[
{\rm with}\;\; \theta_j=\begin{cases} x_j\;\;\;{\rm if}\;\;\; j\in B \\ 0\;\;\;{\rm otherwise}\end{cases} \]
applied to the function
\begin{multline*}
\hat{\varrho}_{k; a_1, \hdots, a_r}(\eta_1, \hdots, \eta_{r-1})=\exp\Bigg(\sum_{j=1}^{r-1}\big(c(k, a_j)-c(k, a_{j+1})\big)\eta_j \Bigg)\\
\times\prod_{\substack{\chi \mod k\\ \chi \neq \chi_0}}F\Bigg(\Bigg|\sum_{j=1}^{r-1}\big(\chi(a_j)-\chi(a_{j-1}) \big)\eta_j \Bigg|, \chi \Bigg)
\end{multline*}
with
\[
F(z, \chi):=\prod_{\substack{\gamma>0\\L(\h+i\gamma, \chi)}}J_0\bigg(\frac{2z}{\sqrt{1/4+\gamma^2}}\bigg)
\]
and 
\[
J_0(z):=\sum_{m=0}^{\infty}\frac{(-1)^m(z/2)^{2m}}{(m!)^2},
\]
 the standard Bessel function of order zero.
\end{thm}
\vskip12pt


%

%

K. Ford and S. Konyagin~\cite{2002.Ford_1} also investigated the Shanks-R\'enyi prime race problem: ostensibly for the races among three competitors: Let $D:=(k, a_1, a_2, a_3)$ where $a_1, a_2, a_3$ are distinct residues modulo $k$ which are coprime to $k$. Suppose for each $\chi \in C_k$ that $B(\chi)$ is a sequence of complex numbers with positive imaginary part (possibly empty, with duplicates allowed), and denote by $\mathcal{B}\,$ the system of $B(\chi)$ for $\chi \in C_k$.  Let $n(\varrho, \chi)$ be the number of occurrences of numbers $\varrho$ in $B(\chi)$.  The system $\mathcal{B}\,$ is called a barrier for $D$ if the following hold:
\begin{enumerate}
\item all numbers in each $B(\chi)$ have real part in $[\beta_2, \beta_3]$, where $1/2<\beta_2<\beta_3 \leq 1$
\item for some $\beta_1$ satisfying $1/2\leq \beta_1<\beta_2$ if we assume that for each $\chi \in C_k$ and $\varrho \in B(\chi)$, $L(s, \chi)$ has a zero of multiplicity $n(\varrho, \chi)$ at $s= \varrho$, and all other zeros of $L(s, \chi)$ in the upper half-plane have real part $\leq \beta_1$, the one of the six ordering of the three functions $\pi_{k, a_i}(x)$ does not occur for large x.
\end{enumerate}
If each sequence $B(\chi)$ is finite, we call $\mathcal{B}$ a finite barrier for $D$ and denote by $|\mathcal{B}|$ the sum of the number of elements of each sequence $B(\chi)$, counted according to multiplicity.
\begin{thm}[\cite{2002.Ford_1} Theorem 1]
For every real number $\tau>0$ and $\sigma>\h$ and for every $D=(k; a_1, a_2, a_3)$, there is a finite barrier for $D$, where each sequence $B(\chi)$ consists of numbers with real part $\leq \sigma$ and imaginary part $> \tau$.  In fact, for most $D$, there is a barrier with $|\mathcal{B}|\leq 3$.
\end{thm}

\vskip12pt

K.Ford and J. Sneed initiated the investigation of biases for products of two primes~\cite{2010.Ford}:
\begin{dfn}
Define $\pi_2(x; k, l)$ to be the number of integers $\leq x$ that are in progression $l \mod k$ and are the product of two (not necessarily distinct) primes, and 
\[
\delta_{\pi_2}(x; k, l_1, l_2):=\pi_2(x; k, l_1)-\pi_2(x; k, l_2)
\]
\end{dfn}
\begin{thm}[\cite{2010.Ford} Theorem 1.1]
Let $a, b$ be distinct elements of $A_k$, where $A_k$ denote the set 
\[
\pi(x; k, a) \sim \frac{x}{\varphi(k)\log x},
\]
then under GRH (Conjecture~\ref{GRH}) and LI (Conjecture~\ref{LI}) for $\chi \mod k$, $\delta_2(k; a, b)$ exists.  Moreover, if $a$ and $b$ are both quadratic residues modulo $q$ or both quadratic non-residues, then $\delta_2(k; a, b)=\h$. (i.e. the race is unbiased)  Otherwise, if $a$ is a quadratic non-residue and $b$ is a quadratic residue, then
\[
1-\delta(k; a, b)<\delta_{\pi_2}(k; a, b)<\h
\]
We can accurately estimate $\delta_2(q; a, b)$ borrowing methods by methods described in~\cite{1994.Rubinstein}.  In particular, we have:
\[
\delta_2(4; 3, 1) \approx 0.10572
\]
\end{thm}

\begin{thm}[\cite{2010.Ford} Theorem 1.2]
Assume GRH (Conjecture~\ref{GRH}) for each $\chi \in C_k$, $L(\h, \chi)\neq 0$ and the zeros of $L(s, \chi)$ are simple.  Then
\[
\frac{\delta_{\pi_2}(x; k, a, b)\log{x}}{\sqrt{x}\log_2 x}=\frac{N_k(b)-N_k(a)}{2\varphi(q)}-\frac{\log x}{\sqrt x}\delta_\pi(x; x, a, b)+\Sigma(x; x, a, b)
\]
where
\[
\frac{1}{Y}\int_1^Y\left|\Sigma(e^y; q, a, b)\right|^2\,dy=o(1) \;\;\;{\rm as}\;\; Y\to \infty,
\]
and $N_k(l)$ is as defined back in~\eqref{asym}.
\end{thm}
\vskip12pt

The most recent developments on the race-problem of Shanks-R\`enyi are due to Y. Lamzouri in his two papers~\cite{201x.Lamzouri_1}~\cite{201x.Lamzouri_2}, where he defines:  
\begin{dfn}
In the notation of Rubinstein and Sarnak~\cite{1994.Rubinstein}, let
\[
\Delta_r(k):=\max_{(a_1, a_2, \hdots, a_r)\in \A_r(k)}\bigg|\delta_{k; a_1, \hdots a_r}-\frac{1}{r!} \bigg|.
\]
\end{dfn}
and he estimates it by:
\begin{thm}[[\cite{201x.Lamzouri_2} Theorem A]
Assume GRH(Conjecture~\ref{GRH}) and LI~(Conjecture~\ref{LI}) for modulo $k$. Let $r\geq 3$ be a fixed integer.\\
If $q$ is large, then
\[
\Delta_r(k) \asymp_r \frac{1}{\log k}.
\]
\end{thm}
He also \emph{redefines} {\bf{unbiased}}:
\begin{dfn}
Let $(a_1, a_2, \hdots, a_r) \in \A_r(k)$, the race $\{q; a_1, \hdots, a_r\}$ is said to be \emph{unbiased} if for every permutation $\sigma$ of the set $\{ 1, 2, \hdots, r \}$ we have
\[
\delta_{q; a_{\sigma(1)}, \hdots, a_{\sigma(r)}}=\delta_{q; a_1, \hdots, a_r}=\frac{1}{r!}.
\]
\end{dfn}

Thus, a race is said to be \emph{biased} if this condition fails to hold, and towards a conjecture made by Rubinstein and Sarnak,

\begin{cnj}[Rubinstein and Sarnak~\cite{1994.Rubinstein}]
When $r\geq 3$, the race $\{q; a_1, \hdots, a_r\}$ is unbiased if and only if $r=3$ and the residue classes $a_1, a_2$ and $a_3$ satisfy the condition 
\begin{equation}
a_2 \equiv a_1\varrho \mod k, \;\;\;\; a_3\equiv a_1 \varrho^2 \mod k,
\end{equation}
for some $\varrho \neq 1$ with $\varrho^3\equiv 1 \mod k$,
\end{cnj}
Lamzouri attacks by,
\begin{thm}[\cite{201x.Lamzouri_2} Theorem B]
Assume GRH~(Conjecture~\ref{GRH}) and LI~(Conjecture~\ref{LI}) modulo $k$. Given $r\geq 3$, there exists a positive number $q_0(r)$ such that for any $k\geq q_0(r)$ there are two $r$-tuples $(a_1, \hdots, a_r), (b_1, \hdots b_r)\in \A_r(k)$, with all the $a_i$'s being squares and all of the $b_i$'s being non-squares modulo $k$, and such that both the races $\{k; a_1, \hdots, a_r\}$ and $\{k; b_1, \hdots, b_r \}$ are biased.
\end{thm}
He also generalizes the definition $\A_r$ made by Rubinstein and Sarnak:
\begin{dfn}
For distinct non-zero integers $a_1, \hdots, a_2$, we define $\mathcal{Q}_{a_1, \hdots, a_r}$ to be the set of positive integers $q$ such that $a_1, \hdots, a_r$ are distinct modulo $q$ and $(q, a_i)=1$ for all $1 \leq i \leq r$.
\end{dfn}
so he could ponder upon
\begin{cnj}[\cite{201x.Lamzouri_2} Conjecture 2]
Let $r\geq 3$ and $a_1, \hdots a_r$ be distinct non-zero integers, then for all positive integers $k\in \mathcal{Q}_{a_1, \hdots, a_2}$ such that $k>2\max(|a_i|^2)$, the race $\{k; a_1, \hdots, a_r \}$ is biased.
\end{cnj}
with
\begin{thm}[\cite{201x.Lamzouri_1} Theorem C]
Assume GRH~(Conjecture~\ref{GRH}) and LI~(Conjecture~\ref{LI}).  Let $r\geq 3$ and $a_1, \hdots, a_r$ be distinct non-zero integers such that one of the following conditions occur:
\begin{enumerate}
\item There exist $1\leq i\neq j \leq r$ such that $a_i+a_j=0$.
\item There exist $1\leq i\neq j \leq r$ such that $a_i/a_j$ is a prime power. 
\end{enumerate}
Then for all but finitely many $k\in\mathcal{Q}_{a_1, \hdots, a_r}$, the race $\{k; a_1, \hdots, a_r\}$ is biased.
\end{thm}
\vskip12pt
Finally, Lamzouri dissects the  measure $\mu_{q; a_1, \hdots, a_r}$ by:
\begin{thm}[\cite{201x.Lamzouri_2} Theorem 1]
Assume GRH (Conjecture~\ref{GRH}) and LI ~(Conjecture~\ref{LI}). For $r\geq 2$ a fixed integer, let $q$ be large and $a_1, \hdots, a_r$ be distinct reduced residues modulo $q$. Then we have 
\[
\mu_{q; a_1, \hdots, a_r}\Big(\|x\|>\lambda \sqrt{\Var(q)} \Big)=(2\pi)^{-r/2}\int_{\|x\|>\lambda}\exp\bigg( -\h \sum_{i=1}^r x_i^2\bigg)\,dx+O_r\bigg(\frac{1}{\log^2 q}\bigg)
\]
for $\lambda$ in the range of $0<\lambda\leq \sqrt{\log_2 q}$.\\
Moreover, there exists an $r$-tuple of distinct reduced classes $(a_1, \hdots,a_r)$ modulo $q$, with $\lambda$ in the range of $1/4<\lambda<3/4$ such that
\[
\bigg|
\mu_{q; a_1, \hdots, a_r}\Big(\|x\|>\lambda\sqrt{\Var(q)} \Big)-(2\pi)^{-r/2}\int_{\| x\|>\lambda}\exp\bigg(-\h\sum_{i=1}^r x_i^2 \bigg)\,dx \bigg| \gg_r \frac{1}{\log^2 q}
\]
\end{thm}
\begin{thm}[\cite{201x.Lamzouri_2} Theorem 2]
Assume GRH(Conjecture~\ref{GRH}) and LI~(Conjecture~\ref{LI}). Fix an integer $r\geq 2$ and a real number $A\geq 1$, $q$ large. For all distinct reduced residues $a_1, \hdots a_r$ modulo $q$, we have
\[
\exp\bigg(-c_{60}(r, A)\frac{V^2}{\varphi(q)\log q}\bigg)\ll \mu_{q; a_1, \hdots, a_r}(\|x\|>V)\ll \exp\bigg(-c_{61}(r, A)\frac{V^2}{\varphi(q)\log q}\bigg)
\]
uniformly within the range $(\varphi(q)\log q)^{1/2}\ll V \leq A\varphi(q)\log q$, where $c_{61}(r, A)> c_{60}(r, A)$ are positive numbers that only depend on $r$ and $A$.
\end{thm}
\begin{thm}[\cite{201x.Lamzouri_2} Theorem 3]
Assume GRH (Conjecture~\ref{GRH}) and LI~(Conjecture~\ref{LI}). For integer $r\geq 2$ and $q$ large, if $V/(\varphi(q)\log q) \to \infty$ and $V/(\varphi(q)\log^2 q) \to 0$ as $q \to \infty$, then for all distinct reduced residues $a_1, \hdots  a_r$ modulo $q$, 
\[\exp\bigg(-c_{63}(r)\frac{V^2}{\varphi(q)\log q}\exp\bigg(c_{65}(r)\frac{V}{\varphi(q)\log q}\bigg)\bigg) \ll \mu_{q; a_1, \hdots, a_r}(\|x\|>V)
\]
and
\[
\mu_{q; a_1, \hdots, a_r}(\|x\|>V) \ll \exp \bigg(-c_{62}(r)\frac{V^2}{\varphi(q)\log q}\exp\bigg(c_{64}(r)\frac{V}{\varphi(q)\log q} \bigg)\bigg)
\]
where $c_{63}(r)>c_{62}(r)$, and $c_{65}(r)>c_{64}(r)$ are positive numbers only depend on $r$.
\end{thm}
\begin{thm}[\cite{201x.Lamzouri_2} Theorem 4]
Assume GRH (Conjecture~\ref{GRH}) and LI~(Conjecture~\ref{LI}). For $q$ large, let $r$ with $2\leq r \leq \varphi(q)-1$ be an integer.  If $V/(\varphi(q)\log^2 q)\to \infty$ as $q \to \infty$, then for all distinct reduced residue classes $a_1, \hdots a_r$ modulo $q$, the tail $\mu_{q; a_1, \hdots a_r}(|x|_\infty > V)$ equals
\[
\exp\Bigg(-e^{L(q)}\sqrt{\frac{2(\varphi(q)-1)V}{\pi}} \exp\bigg(\sqrt{L(q)^2+\frac{2\pi V}{\varphi(q)-1}}\bigg)\bigg(1+O\bigg(\Big(\frac{\varphi(q)\log^2(q)}{V} \Big)^{1/4} \bigg) \bigg)\Bigg),
\]
where 
\[
L(q)=\frac{\varphi(q)}{\varphi(q)-1}\bigg(\log q -\sum_{p|q}\frac{\log p}{p-1} \bigg) + A_0 -\log \pi,
\]
and 
\[
A_0:=\int_0 ^1\frac{\log I_0(t)}{t^2}\,dt +\int_1 ^\infty \frac{\log I_0(t)-t}{t^2}\,dt +1, 
\]
with $\displaystyle I_0(t)=\sum_{n=0}^\infty \frac{(t/2)^{2n}}{(n!)^2}$ being the modified Bessel function of order zero.
\end{thm}
giving a conditional bound of the tails to the measure $\mu_{q; a_1, \hdots, a_r}$, fully generalizing the work done by Montgomery~\cite{1980.Montgomery} on $\mu_{1; 1}$.

\end{document}